\documentclass[11pt,a4paper]{article}
\newcommand{\vect}[1]{\boldsymbol{#1}}
\newcommand{\vectt}[1]{\boldsymbol{\mathbf{#1}}}

\newcommand{\dx}{\mathrm{d}x}

\newcommand{\supp}{\mathrm{supp}}
\newcommand{\E}{\mathrm{e}}
\newcommand{\fdx}{\frac{\mathrm{d}}{\dx}}
\newcommand{\lt}{{L^2(\Omega)}}

\newcommand{\hoo}{{H^1_0(\Omega)}}
\newcommand{\hp}{{hp}}

\newcommand{\Uh}{U_{h,\vect{p}}}
\newcommand{\Ph}{\Psi_{h,\vect{p}}}

\newcommand{\spa}{\mathrm{span}}

\def\XXint#1#2#3{{\setbox0=\hbox{$#1{#2#3}{\int}$}
\vcenter{\hbox{$#2#3$}}\kern-.5\wd0}}

\usepackage{color}
\definecolor{deepblue}{rgb}{0,0,0.5}
\definecolor{deepred}{rgb}{0.6,0,0}
\definecolor{deepgreen}{rgb}{0,0.5,0}

\usepackage{amsmath,amssymb}
\usepackage{bm}
\usepackage{amsthm}
\usepackage{anyfontsize}
\usepackage{url}
\usepackage{graphicx}
\usepackage{array}
\usepackage{varwidth}
\usepackage{geometry}
\usepackage{ esint }
\definecolor{ao(english)}{rgb}{0.0, 0.5, 0.0}
\usepackage{listings}
\usepackage{fancyhdr}
\usepackage{enumitem}
\usepackage[
    backend=bibtex,
   bibstyle=ieee,
    citestyle=numeric-comp,
    natbib=true,
    sortlocale=en_GB,
    sorting = nyt,
    url=true,
    doi=true,
    arxiv=true,
    eprint=true
]{biblatex}
\DeclareFieldFormat{eprint:arXiv}{%
  arXiv preprint: \href{https://arxiv.org/abs/#1}{#1 \texttt{[\printfield{eprintclass}]}}%
}

\AtBeginBibliography{\small}

\usepackage[table]{xcolor}
\usepackage{hhline}

\usepackage[font=small]{caption}
\usepackage{mathtools}
\usepackage{subfig}
\usepackage{appendix}
\usepackage{standalone}
\usepackage{tikz}
\usetikzlibrary{trees}
\usepackage{hyperref}
\usepackage{booktabs}

\hypersetup{
    colorlinks,
    linkcolor={blue!50!black},
    citecolor={blue!50!black},
    urlcolor={blue!80!black}
}

\RequirePackage{algorithm}[0.1]

\RequirePackage[capitalize,nameinlink]{cleveref}[0.19]
\crefformat{equation}{\textup{#2(#1)#3}}
\crefrangeformat{equation}{\textup{#3(#1)#4--#5(#2)#6}}
\crefmultiformat{equation}{\textup{#2(#1)#3}}{ and \textup{#2(#1)#3}}
{, \textup{#2(#1)#3}}{, and \textup{#2(#1)#3}}
\crefrangemultiformat{equation}{\textup{#3(#1)#4--#5(#2)#6}}%
{ and \textup{#3(#1)#4--#5(#2)#6}}{, \textup{#3(#1)#4--#5(#2)#6}}{, and \textup{#3(#1)#4--#5(#2)#6}}

\Crefformat{equation}{#2Equation~\textup{(#1)}#3}
\Crefrangeformat{equation}{Equations~\textup{#3(#1)#4--#5(#2)#6}}
\Crefmultiformat{equation}{Equations~\textup{#2(#1)#3}}{ and \textup{#2(#1)#3}}
{, \textup{#2(#1)#3}}{, and \textup{#2(#1)#3}}
\Crefrangemultiformat{equation}{Equations~\textup{#3(#1)#4--#5(#2)#6}}%
{ and \textup{#3(#1)#4--#5(#2)#6}}{, \textup{#3(#1)#4--#5(#2)#6}}{, and \textup{#3(#1)#4--#5(#2)#6}}

\usepackage{algpseudocode}
\algnewcommand{\Initialize}[1]{
  \State \textbf{Initialize:}
  \Statex \hspace*{\algorithmicindent}\parbox[t]{.8\linewidth}{\raggedright #1}
}
\algnewcommand{\Indent}[2]{
  \State {#1}
  \vspace{-2mm}
  \Statex \hspace*{\algorithmicindent}\parbox[t]{.9\linewidth}{\raggedright #2}
}

\usepackage{todonotes}

\newtheorem{theorem}{Theorem}[section]

\newtheorem{lemma}[theorem]{Lemma}
\newtheorem{definition}[theorem]{Definition}

\newtheorem{remark}[theorem]{Remark}
\numberwithin{equation}{section}

\usepackage{todonotes}

\title{Hierarchical proximal Galerkin: a fast $hp$-FEM solver for variational problems with pointwise inequality constraints}

\author{Ioannis P.~A.~Papadopoulos\thanks{\small Mathematical Institute, University of Oxford, UK,  \tt{ioannis.papadopoulos@maths.ox.ac.uk}.}}

\begin{document}

\maketitle
\date
\thispagestyle{empty}
\pagestyle{fancy}
\lhead{I.~P.~A.~Papadopoulos}

\begin{abstract}
We leverage the proximal Galerkin algorithm (Keith and Surowiec, Foundations of Computational Mathematics, 2024), a recently introduced mesh-independent algorithm, to obtain a high-order finite element solver for variational problems, posed on tensor-product domains, with pointwise inequality constraints. This is achieved by discretizing the saddle point systems, arising from the latent variable proximal point method, with the hierarchical $p$-finite element basis. This results in discretized sparse Newton systems that admit a simple and effective block preconditioner. The solver can handle both obstacle-type, $u \leq \varphi$, and gradient-type, $|\nabla u| \leq \varphi$, constraints.  We apply the resulting algorithm to solve obstacle problems with $hp$-adaptivity, a three-dimensional obstacle problem, a gradient-type constrained problem, and the thermoforming problem, an example of an obstacle-type quasi-variational inequality. We observe $hp$-robustness in the number of Newton iterations and only mild growth in the number of inner Krylov iterations to solve the Newton systems. Crucially we also provide wall-clock timings that are faster than low-order discretization counterparts.
\end{abstract}

\section{Introduction}
\label{sec:introduction}

Minimizing the Dirichlet energy with an obstacle- or gradient-type pointwise constraint is known as the \emph{obstacle problem} or \emph{generalized elastic-plastic torsion problem}, respectively. The obstacle problem models the deflection of an elastic membrane in contact with an obstacle \cite{fichera1964, stampacchia1964} and due to its foundational nature, has wide-ranging applications in various fields including optimal control \cite{christof2021,kinderlehrer2000,rodrigues1987}, topology optimization \cite{Bendsoe2004,Borrvall2003,Papadopoulos2021a}, and elasticity theory \cite{hlavacek1988, fichera1964} among many others. In turn, many PDE-constrained optimization problems are thoroughly reliant on an efficient obstacle problem solver. Gradient-type constraints arise in quasi-variational inequalities that model  sandpile growth, river networks, semiconductors, strain-limited elastic material, and number of other applications cf.~\cite{arndt2021,white1983,TWTing1969,Brezis1971a,Prigozhin1996,Prigozhin1996a,Antil2022,
kunze2000,evans1979,santos2002,ting1977, caffarelli1979, rodrigues2019}.

In this paper we discretize the latent variable proximal point method (LVPP) of Keith and Surowiec \cite{keith2023} with the hierarchical $p$-finite element method ($p$-FEM) to construct the hierarchical proximal Galerkin algorithm (hpG), a solver that enjoys the following properties:
\begin{enumerate}[label=(\roman*)]
\itemsep=-2pt
\item $h$- and $p$-independent number of nonlinear iterations until convergence;
\item arbitrarily fast convergence to solutions;%
\item high-order discretizations (see \Cref{sec:examples:qvi} for a discretization with polynomials up to total degree 164 $(p=82)$ on each quadrilateral element of a $4 \times 4$ mesh);
\item competitive solve times, when the problem exhibits tensor-product structure, with lower order methods to reach a prescribed error;
\item a block-diagonal preconditioner resulting in mild polylogarithmic growth in the inner Krylov method iterations as $h\to0$ and $p\to\infty$;
\item fast sparse factorizations for the mass, stiffness, and preconditioner matrices;
\item fast quasi-optimal complexity quadrature via fast transforms and the discrete cosine transform (DCT).
\end{enumerate}

\subsection{Motivation}
For nontrivial constraints, the regularity of a solution $u$ is capped at $u \in H^{s}(\Omega)$, for $s<5/2$, where $H^s(\Omega)$ denotes the usual Sobolev spaces $W^{s,2}(\Omega)$ \cite{Adams2003}. This is true even if the problem is posed on a smooth domain  $\Omega$ with smooth data and a smooth obstacle \cite[Th.~3]{brezis1971}\footnote{For a constant gradient-type constraint, the fact that $u \in H^s(\Omega)$, for $s<5/2$, follows by rewriting the gradient-type constrained problem as an obstacle problem with a Lipschitz continuous obstacle \cite[Th.~2.2]{ting1977}. The result in \cite[Th.~3]{brezis1971} then applies.}. As such low-order discretizations are often favoured due to their simplicity, perceived reduced computational cost due to sparsity, and efficiency in convergence due to low solution regularity. For the obstacle problem, there exist efficient low-order optimal complexity (but mildly mesh dependent) multigrid (MG) techniques \cite{graser2009,brandt1983,hackbusch1983,hoppe1987,hoppe1994,kornhuber1994, kornhuber1996, bueler2024, wohlmuth2003}.  There is a strong correspondence between these MG solvers and active-set strategies \cite[Sec.~6]{graser2009} which require that nodal feasibility in the discretization implies pointwise feasibility. In turn this requirement typically restricts a user to the lowest-order continuous FEM basis\footnote{A recent study has shown success with active-set strategies and a carefully chosen quadratic or cubic FEM basis \cite{kirby2024}.}. A reader may now arrive at the following question:

\vspace{-3mm}
\begin{center}
\textit{Why develop a high-order discretization for pointwise-constrained problems?}
\end{center}
\vspace{-3mm}

\smallskip

\noindent In this paper, we argue that, per degree of freedom (dof), a high-order discretization may induce an error that is orders of magnitude smaller than the low-order counterpart -- particularly in the $H^1$-norm. Of course, fewer dofs does not imply faster solve times. Matrix assembly overhead may significantly increase as sparsity is often greatly reduced, resulting in a substantial increase in the number of nonzero entries that must be evaluated. Moreover, each individual nonzero entry becomes more expensive to compute. The loss of sparsity also impacts the effectiveness of sparse linear solvers further reducing solver speed.  However, in recent years, there have been a number of developments that allow for fast expansions and evaluations of the high-order orthogonal polynomials that form the hierarchical $p$-FEM basis which, incidently, when used to discretize the LVPP subproblems \cite{keith2023}, induces discretized linear systems where the loss of sparsity is heavily mitigated \cite{townsend2018, knook2024quasi,papadopoulos2024b}. By leveraging these developments, we will present concrete counterexamples where accurate solutions are cheaper to obtain with high-order discretizations when the problem exhibits tensor-product structure.

Secondly, despite the cap on the regularity of the solution, convergence rates faster than $O((h/p)^{3/2})$ in the $H^1$-norm are achievable if the mesh aligns with the regions in the solution where the loss of regularity occurs. In our context, these regions are where the pointwise constraints transition from being active to inactive. These are not known beforehand. However, one may develop estimators to locate these regions and, particularly in the case of gradient-type constraints, the location of the transitions may sometimes be inferred from the problem data. In \Cref{sec:examples:1d} and \Cref{sec:examples:gradient} we present two examples where we obtain convergence rates faster than $3/2$ in the regimes of $h$ and $p$ we consider. Note that a discretization with degree $p=1$ will only ever converge at a maximum rate of $O(h)$.

Aside from solver speed and convergence rates, high-order discretizations often offer other benefits. They minimize numerical artefacts in fluid dynamics and elasticity \cite{Karniadakis2005}, e.g.~locking in linear elasticity does not occur if $p\geq 2d$, $d \in \{2,3\}$ \cite{Ainsworth2022}. As pointwise constraints naturally occur in models for such settings (e.g.~pointwise constraints on the density field in the topology optimization of linear elasticity \cite{Bendsoe2004,papadopoulos2024c} and fluids \cite{Borrvall2003,papadopoulos2022a,papadopoulos2022d}) then utilizing a solver that is then confined to the lowest-order (continuous) FEM discretization is limiting. High-order methods are also amenable to parallelization making them suitable for modern computing architectures cf.~\cite{knook2024quasi,andrej2024}.

The proximal Galerkin (pG) algorithm arises after a Galerkin discretization of LVPP which was first introduced in the context of the obstacle problem by Keith and Surowiec \cite{keith2023} and generalized to a number of other problems in a recent review \cite{dokken2024}. LVPP is formulated on the infinite-dimensional level and not tied to any particular choice of discretization. The main component of the algorithm is the repeated solve of a coupled nonlinear system of PDEs. Both theoretical and numerical investigations indicate that the number of iterations of the pG solver is mesh ($h$) and degree ($p$) independent \cite{keith2025}. This, among other reasons, is a motivating factor for choosing a variation of the pG solver in this work. A second reason is that after a Newton linearization of the pG nonlinear system, one arrives at a saddle point problem with one nontrivial term. Hence, with a sparsity promoting discretization, the density of a high-order method is localized to the nontrivial term which in turn is effectively handled by a simple block preconditioning strategy.

In principle, a high-order FEM discretization of a path-following penalty method can be used to enforce the pointwise constraints. In exact arithmetic, this would lead to mesh independent iteration counts \cite{hintermuller2009,hintermuller2004}. However, in order to extract optimal convergence rates for the obstacle problem when $p > 1$, one must scale the penalty parameter as $O(h^{-3})$ which also scales the condition number of the arising discretized linear systems by the same factor \cite[Sec.~3.1]{gustafsson2017}. Ill-conditioning may quickly cause mesh dependent effects and a potential blowup in the required number of nonlinear iterations. Moreover, it is unclear how to handle the loss of sparsity as $p \to \infty$. 

We are not the first to discuss high-order discretizations for the obstacle problem \cite{porwal2024}. Keith and Surowiec, in the seminal pG paper, solve a 2D obstacle problem with a mesh consisting of five cells and a degree $p=12$ basis \cite[Fig.~11]{keith2023}. They were capped at $p=12$ by their choice of FEM software and their chosen basis induces increasingly dense stiffness and mass matrices as $p\to \infty$. They solved the linear systems via a direct solver. A more recent paper on pG includes an obstacle problem example with $p=48$ although the discretization is a spectral method limited to a single-cell disk mesh \cite[Sec.~3.1]{dokken2024}. An iterative solver was used but resulted in faster growing iteration counts than what we observe in our methods. Other works include \cite{krebs2007,  gwinner2013} as well as a sequence of successful papers by Banz, Schr\"oder, and coauthors \cite{banz2015, banz2014, banz2020, banz2021}. In these works the choice of FEM basis also leads to increasingly denser stiffness and mass matrices as $p \to \infty$, the outer nonlinear solver is mesh dependent, i.e.~the number of nonlinear iterations grow as $h \to 0$ and $p \to \infty$, and the implementation relies on direct solvers. However, in \cite{banz2015}, Banz and Schr\"oder successfully consider examples with $p=30$ and implement a fully $h$- and $p$-adaptive scheme in two dimensions, complete with a posteriori estimators, leading to highly accurate approximations. Also of significant note is a $hp$-FEM discretization to a problem involving gradient-type constraints in \cite{bammer2024}, where Bammer, Banz, and Schr\"oder considered a discretization with $p=50$.

\subsection{Contributions}

The main contribution of this paper is the investigation of a fast sparsity-preserving $hp$-FEM discretization of the subproblems of the LVPP algorithm applied to the obstacle and generalized elastic-plastic torsion problems on tensor-product domains, with a view of providing a bedrock for its potential in applications. We develop the first preconditioners for the Newton systems arising in a pG method and observe them to be exceptionally effective. These appear to also be the first preconditioners for high-order discretizations of obstacle and gradient-type constraints. Classical works follow discretize-then-optimize approaches and are therefore confined to the lowest order discretization \cite{pearson2017,rees2010}.

The quadrature for the assembly of the matrices at each Newton step can be confined to a discontinuous basis consisting of Legendre polynomials which admit fast evaluation and expansion transforms. We then leverage an atypical but fast quadrature scheme. We also see that the method may adopt error estimators already developed in the literature to guide $hp$-refinement. Altogether this allows us to consider examples with the highest discretization degree that have been reported in the literature to date. Moreover, since the novel hpG solver retains sparsity, requires little computational overhead, and exhibits  $hp$-independent number of nonlinear iterations, we are able to provide the first reported competitive wall-clock timing comparisons with other solvers featuring low-order discretizations. In fact in \Cref{sec:examples}, we observe up to a 24 and 100 times speed up, per linear solve, to reach the same error for an obstacle problem and generalized elastic-plastic torsion problem, respectively (cf.~\Cref{sec:examples:2d:1,sec:examples:gradient}).

We begin by introducing the Dirichlet minimization problem with obstacle- and gradient-type constraints, together with the pG algorithm in \Cref{sec:pg}, followed the construction of the hierarchical $p$-FEM basis in \Cref{sec:fem}. Next we provide a detailed investigation of the discretized Newton linearization of the pG subsystems in \Cref{sec:newton} and the design of a preconditioning strategy in \Cref{sec:preconditioning}. We outline four algorithmic techniques that significantly speed up the algorithm when $p \gg 1$ in \Cref{sec:implementation}. We consider five numerical examples in \Cref{sec:examples}. \Cref{sec:examples:1d} tackles a one-dimensional obstacle problem with an oscillatory forcing-term where  $hp$-adaptivity is fully explored.  \Cref{sec:examples:2d:1} handles a two-dimensional obstacle problem where it is clear that high-order discretizations give smaller errors both with respect to the dofs \emph{and} wall-clock time (see \cref{fig:oscillatory-obstacle}). Next, in \Cref{sec:examples:gradient}, we investigate the effectiveness of the solver applied to a gradient-type constrained problem. The penultimate problem in \Cref{sec:examples:qvi} computes the approximate solution of a two-dimensional thermoforming problem, an obstacle-type quasi-variational inequality, with a discretization consisting of 16 cells and multivariate polynomials up to total degree 164 on each element. We conclude the examples with the solve of a three-dimensional obstacle problem in \Cref{sec:examples:3d}. We compare various iterative solvers for the underlying linear system that arise. Finally we give our conclusions in \Cref{sec:conclusions}.%

\section[Setup]{Problem setup and latent variable proximal point}
\label{sec:pg}
Consider an open, convex, and bounded Lipschitz domain $\Omega \subset \mathbb{R}^d$, $d \in \{1,2\}$. Let $\langle \cdot, \cdot \rangle_{X^*,X}$ denote the (topological) duality pairing between a Banach space $X$ and its dual space $X^*$. The inner product of a Hilbert space $H$ is denoted by $(\cdot, \cdot)_H$.  Let $W^{s,p}(\Omega)$, $s > 0$, $p \geq 1$, denote the usual Sobolev spaces and let $L^p(\Omega)$, $p\geq 1$, denote the Lebesgue spaces \cite{Adams2003}. Consider the Hilbert spaces $H^s(\Omega) \coloneqq W^{s,2}(\Omega)$, $s > 0$ and denote by $H^1_0(\Omega) \coloneqq \{v \in H^1(\Omega) : v|_{\partial \Omega}= 0\}$ where $|_{\partial \Omega} : W^{1,p}(\Omega) \to W^{1-\frac{1}{p}, p}(\partial \Omega)$ denotes the standard boundary trace operator \cite{Gagliardo1957}. Consider the Dirichlet energy functional
\begin{align}
J(v) \coloneqq \frac{1}{2} \| \nabla v\|^2_{L^2(\Omega)} - (f, v)_{L^2(\Omega)}.
\label{eq:vi1}
\end{align}
We seek the minimizer $u \in H^1_0(\Omega)$ of the Dirichlet energy functional $J$, with datum $f \in L^2(\Omega)$, while satisfying the obstacle or gradient-type constraint with $\varphi \in H^{1}(\Omega)$. In other words we consider one of the obstacle- or gradient-type constrained minimization problems:
\begin{subequations}
\begin{align}
\min_{u \in K} J(u) \;\; \text{where} \;\; K &= \{ v \in H^1_0(\Omega) : v \leq \varphi \; \text{a.e.~in} \; \Omega \}&&\text{(obstacle-type),} \label{eq:ad:obstacle}\\
\min_{u \in K} J(u) \;\; \text{where} \;\; K &=\{ v \in H^1_0(\Omega) : |\nabla v| \leq \varphi \; \text{a.e.~in} \; \Omega \} && \text{(gradient-type)}.\label{eq:ad:gradient}
\end{align}
\end{subequations}

\subsection{LVPP: obstacle problem}
LVPP approximates the solution of the obstacle problem \cref{eq:ad:obstacle}, on the infinite-dimensional level, by solving a series of nonlinear mixed saddle point systems where a parameter $\alpha \in \mathbb{R}_+$ and a source term is updated. More specifically, given the datum $f \in L^\infty(\Omega)$, the obstacle $\varphi \in H^{1}(\Omega)$, $\varphi|_{\partial \Omega} \geq 0$, the parameter $\alpha_{k} > 0$, $k \in \mathbb{N}$, a previous latent variable iterate $\psi_{k-1} \in L^\infty(\Omega)$ and setting $\psi_0 \equiv 0$, the pG subproblem seeks $(u_{k}, \psi_{k}) \in H^1_0(\Omega) \times L^\infty(\Omega)$ that satisfies, for all $(v, \zeta) \in \hoo \times L^\infty(\Omega)$ \cite[Alg.~3]{keith2023},
\begin{align}
\begin{split}
\alpha_{k} (\nabla u_{k}, \nabla v)_\lt + (\psi_k, v)_\lt &= \alpha_k(f, v)_\lt + (\psi_{k-1}, v)_\lt,\\
(u_k, \zeta)_\lt + (\E^{-\psi_k}, \zeta)_\lt &= (\varphi, \zeta)_\lt.
\end{split} \label{eq:pG}
\end{align}
It was shown that as long as $\sum_{k=1}^N \alpha_k \to \infty$ as $N \to \infty$ then $u_k \to u^*$ strongly in $H^1(\Omega)$ and $\lambda_k \coloneqq (\psi_{k-1} - \psi_k)/\alpha_k \to \lambda^*$ strongly in $H^{-1}(\Omega)$ where $u_*$ is the solution of the obstacle problem and $\lambda^* \coloneqq -\Delta u^* - f$ is its associated Lagrange multiplier \cite[Th.~4.13]{keith2023}. By driving $\alpha_k$ larger at a faster-than-geometric rate (e.g.~$\alpha_k = k(k+1)(k+2)\cdots(k+m)$ for some $m \in \mathbb{N}$), one can converge to the solution of the obstacle problem with superlinear convergence \cite[Cor.~A.12]{keith2023}. Moreover, unlike other mesh independent solvers that rely solely on a penalty term, cf.~\cite{hintermuller2006,hintermuller2009}, the solver still converges to the solution even if $\alpha_k$ is kept at a fixed value, albeit at a sublinear rate \cite[Cor.~A.12]{keith2023}. The faster $\alpha_k$ increases, the quicker the convergence to the solution of the original obstacle problem but at the cost of a harder nonlinear problem at each iteration.

\subsection{LVPP: gradient-type constraint}
Recently LVPP was adapted to handle gradient-type constraints \cite[Sec.~4.1]{dokken2024} although the convergence properties and choices of stable discretizations have not yet been theoretically justified. Nevertheless, experimentally, we observe very similar behaviours as for the obstacle problem. As before, suppose we are given the datum $f \in L^\infty(\Omega)$, the constraint function $\varphi \in L^\infty(\Omega) \cap H^{1}(\Omega)  $ where $\mathrm{ess \, inf} \, \varphi \geq c$ for some $c>0$, the parameter $\alpha_{k} > 0$, $k \in \mathbb{N}$, and a previous latent variable iterate ${\psi}_{k-1} \in L^\infty(\Omega)^d$ where if $k=1$ then ${\psi}_0 \equiv 0$. Then, the pG gradient-type subproblem for approximating the solution of \cref{eq:ad:gradient} is to find $(u_{k}, {\psi}_{k}) \in H^1_0(\Omega) \times L^\infty(\Omega)^d$ that satisfies, for all $(v, {\zeta}) \in \hoo \times L^\infty(\Omega)^d$ \cite[Eq.~(4.3)]{dokken2024},
\begin{align}
\begin{split}
\alpha_{k} (\nabla u_{k}, \nabla v)_\lt + ({\psi}_k, \nabla v)_\lt &= \alpha_k(f, v)_\lt + ({\psi}_{k-1}, \nabla v)_\lt,\\
(\nabla u_k, {\zeta})_\lt &= ({\varphi {\psi}_k}{({1+|{\psi}_k|^2})^{-1/2}},{\zeta})_\lt.
\end{split} \label{eq:pG:gradient}
\end{align}

\section{Hierarchical $p$-FEM basis}
\label{sec:fem}
The hierarchical $p$-FEM basis was pioneered in the 1980s and 1990s by Babu{\v{s}}ka, Szab{\'o} and coauthors \cite{Babuska1981a, Babuska1981b, Babuska1991, babuvska1983lecture,szabo2011introduction}. The basis contains very high-order piecewise polynomials whilst retaining sparsity in the stiffness (weak Laplacian) and mass matrices. More specifically, the basis functions consist of shape and external functions, otherwise known as bubble and hat functions. The shape functions are translated and scaled weighted orthogonal polynomials whose support are contained within a single cell in the mesh. Hence, shape functions tested against shape functions supported on another cell evaluate to zero. Moreover, thanks to the orthogonality properties, testing many of the basis functions with other basis functions supported on the same cell also evaluate to zero. For more details we refer the reader to \cite{Schwab1998,Karniadakis2005,knook2024quasi}.

In this work we consider both a continuous $H^1$-conforming hierarchical $p$-FEM space, which we denote $\Uh$ to discretize $u$ and a discontinuous $L^2$-conforming space denoted by $\Ph$ to discretize $\psi$. We first define these spaces on a one-dimensional domain $\Omega = (a,b)$.

Consider the one-dimensional mesh $\mathcal{T}_h = \{x_i\}_{i=1}^m$ where $a = x_0 < x_1 < \cdots < x_m =b$ and $h = \max_{i} |x_{i+1} - x_i|$.  Let $K_i=[x_i,x_{i+1}]$ and consider the set
\begin{align}
\Psi_{K_i,p} \coloneqq \left\{ P_n(y) : y = \frac{2x-x_i-x_{i+1}}{x_{i+1}-x_i}, \; 0 \leq n \leq p \right\},
\end{align}
where $P_n(x)$ denotes the Legendre polynomial of degree $n$ supported on the interval $[-1,1]$ \cite[Sec.~18.3]{NIST:DLMF}. For $x \in \mathbb{R} \backslash [-1,1]$ we define $P_n(x) = 0$. The Legendre polynomials are orthogonal with respect to the $L^2$-inner product, i.e.~$(P_n, P_m)_{L^2(-1,1)} = \frac{2 \delta_{nm}}{2n+1}$, where $\delta_{nm}$ is the Kronecker delta. 
\begin{definition}[$L^2$-conforming hierarchical $p$-FEM space $\Ph$]
Consider the one-dimensional mesh $\mathcal{T}_h = \cup_{i=1}^m K_i$. We define the $L^2$-conforming hierarchical $p$-FEM space $\Ph$, $\vect{p} = (p_1, \dots, p_m)$ as the union of the cellwise bases:
\begin{align}
\Ph \coloneqq \spa \; \bigcup_{i=1}^m \Psi_{K_i,p_i}.
\end{align}
In other words $\Ph$ consists of scaled and translated Legendre polynomials on cell $K_i$ up to degree $p_i$ where each basis function has a support that is contained in a single cell on the mesh. 
\end{definition}
\begin{lemma}[Diagonal mass matrix]
\label{lem:diagonal-mass-matrix}
The mass matrix of the discontinuous finite element space $\Ph$ is diagonal, i.e.~$(\zeta_i, \zeta_j)_\lt = \frac{|K_i| \delta_{ij}}{2i + 1} $ for all basis functions $\zeta_i, \zeta_j \in \Ph$ and $K_i \in \mathcal{T}_h$ is the cell such that $\supp(\zeta_i) \subseteq K_i$.\end{lemma}
The result in \cref{lem:diagonal-mass-matrix} follows as a direct consequence of the fact that each basis function is only supported on one cell on the mesh as well as the orthogonality of the basis functions supported on the same cell with respect to the $L^2$-norm.

The $H^1$-conforming space $\Uh$ consists of the union of the set of standard piecewise linear hat functions with weighted Jacobi polynomials.  Given two adjacent cells $K_{i-1}$, $K_i$, let $H_i(x)$ denote the usual continuous piecewise linear hat function supported on  $K_{i-1} \cup K_i$ \cite[Ch.~0.4]{Brenner2008}.
Moreover, consider the polynomials
\begin{align}
W_n(x) \coloneqq \frac{(1-x^2) P_n^{(1,1)}(x)}{2(n+1)},
\label{eq:Wn}
\end{align}
where $P_n^{(1,1)}(x)$ are the Jacobi polynomials orthogonal with respect to the weight $1-x^2$ on the interval $[-1,1]$ \cite[Sec.~18.3]{NIST:DLMF}. The normalization is chosen such that $\fdx W_n(x) = -P_{n+1}(x)$ for $n \geq 0$. We define the set of shape functions for the cell $K_i$ as
\begin{align}
B_{K_i,p} \coloneqq \left\{ W_n(y) : y = \frac{2x-x_i-x_{i+1}}{x_{i+1}-x_i}, \; 0 \leq n \leq p-2 \right\}.
\end{align}
\begin{definition}[$H^1$-conforming hierarchical $p$-FEM space $\Uh$]
Consider the one-dimensional mesh $\mathcal{T}_h = \cup_{i=1}^m K_i$. We define the $H^1$-conforming hierarchical $p$-FEM space $\Uh$, $\vect{p} = (p_1,\dots,p_m)$, as:
\begin{align}
\Uh \coloneqq \spa \; \left\{\{H_i\}_{i=1}^{m+1} \cup \bigcup_{i=1}^m B_{K_i,p_i}\right\}.
\end{align}
We use $U_{0,h,\vect{p}}$ to denote the space
\begin{align}
U_{0,h,\vect{p}} \coloneqq \{ v_\hp \in U_{h,\vect{p}} : v_\hp|_{\partial \Omega} = 0\}.
\end{align}
\end{definition}

 We take the tensor product space for two and three dimensions.
We refer to the degree of the pre-tensor one-dimensional polynomials as the \emph{partial degree} and the degree of the resulting multivariate polynomial as the \emph{total degree}. With a slight abuse of notation, we will often compactify the subscripts such that $U_\hp$ and $\Psi_\hp$ denote the the $H^1$-conforming and discontinuous FEM spaces, respectively, in 1D, 2D, and 3D.

We take the opportunity to introduce a discontinuous basis consisting of Legendre spectral Galerkin polynomials \cite{shen1994}. This basis will later be used in the stabilization term for gradient-type constraints as well as to construct a preconditioner for the obstacle problem. Consider a one-dimensional reference cell $K=[-1,1]$. We define the spectral Galerkin polynomials, for $n \in \mathbb{N}_0$:
\begin{align}
Y_n(x) \coloneqq P_n(x) - P_{n+2}(x).
\end{align}
Note that $Y_n(\pm 1) = 0$ for all $n \in \mathbb{N}_0$.
In fact each polynomial $Y_n$ is a rescaling of the shape function $W_n$ defined in \cref{eq:Wn}, i.e.~$Y_n(x) = c_nW_n(x)$ for some $c_n \in \mathbb{R}$. Nevertheless, we found that the numerical convergence of the linear system iterative solvers proposed later were often better when utilizing the rescaling rather than the shape functions directly. Hence we choose to keep the additional notation $Y_n$ for clarity on where a rescaling is used or not.
We denote the FEM basis consisting of these basis functions on each cell as $\Phi_{h, p}$. As before, the two- and three-dimensional basis is constructed via a tensor-product.

\section[Proximal Galerkin]{Proximal Galerkin, Newton systems \& preconditioning}
\label{sec:newton}

We employ a mixed $H^1 \times L^2$-conforming hierarchical finite element discretization of the LVPP saddle-point problems for $(u,\psi)$ in \cref{eq:pG,eq:pG:gradient} where, in particular, the discretization is continuous for $u$ ($u_\hp \in U_\hp$) and discontinuous for $\psi$ ($\psi_\hp \in \Psi_\hp$ in \cref{eq:pG} and $\psi_\hp \in \Psi_\hp^d$ in \cref{eq:pG:gradient}). We refer to the discretized nonlinear systems as the hpG subproblems. Below we first discuss the discretization of the obstacle-type LVPP subproblem \cref{eq:pG} and detail the differences for the gradient-type LVPP subproblems in \Cref{sec:newton:gradient}.

\subsection{Obstacle problem}
\label{sec:newton:obstacle}
For the obstacle-type LVPP subproblem \cref{eq:pG}, we discretize $u$ with (partial) degree $\vect{p}$ and $\psi$ with (partial) degree $\vect{p}-2$ for $\vect{p}\geq 2$, where $\pm 1$ and $\geq$ are to be understood entrywise. With a slight abuse of notation, we denote the finite element functions approximating $u$ and  $\psi$ 
by $u_\hp$ and $\psi_\hp$, respectively. 

The choice of degree pairing for $u$ and $\psi$ is motivated by the inf-sup stability result for $p$-continuous and $(p-2)$-discontinuous FEM on shape-regular sequences of affine meshes in \cite[Lem.~B.3]{keith2023}.  As we will show, by picking a discontinuous discretization for the $\psi$, the linearization of the nonlinearity in \cref{eq:pG} is approximated by a block-diagonal matrix. Otherwise, for a continuous discretization, as $p$ grows, the FEM matrix would become dense. We remark that equal-order continuous discretizations for $u$ and $\psi$ are also inf-sup stable, albeit on quasi-uniform meshes \cite[Remark 5.2]{keith2023}. Further results on stable FEM pairs include \cite{ern2026,fu2024,keith2025}.

After a Newton linearization of the FEM discretization of \cref{eq:pG}, a routine derivation reveals that one must (approximately) solve the following symmetric saddle-point linear system at each Newton step:
\begin{align}
G \coloneqq
\begin{pmatrix}
A_\alpha & B\\
B^\top & -D_\psi - E_\beta 
\end{pmatrix}
\begin{pmatrix}
\vectt{\delta}_u\\
\vectt{\delta}_\psi
\end{pmatrix}
= 
\begin{pmatrix}
\vectt{b}_u\\
\vectt{b}_\psi
\end{pmatrix},
\label{eq:pG-matrix}
\end{align}
where $\vectt{b}_u$ and $\vectt{b}_\psi$ denote the $u$ and $\psi$ components of the residual and $\vectt{\delta}_u$ and $\vectt{\delta}_\psi$ are the coefficient vectors for the Newton updates of $u_\hp$ and $\psi_\hp$, respectively. $A_\alpha \coloneqq \alpha A$ is the $\alpha$-scaled stiffness matrix for $U_\hp$, $B$ is the Gram matrix between $U_{0,\hp}$ and $\Psi_\hp$, $D_\psi$ is the discretization of the linearization of the nonlinear term, and $E_\beta$ is a symmetric positive-definite stabilization matrix such that
\begin{equation}
\begin{aligned}
A_{ij} &= (\nabla {v}_i, \nabla {v}_j)_\lt, \quad  &[D_\psi]_{ij} &= (\zeta_i, \E^{-\psi_\hp} \zeta_j)_\lt,\\
B_{ij} &=  ({v}_i, {\zeta}_j)_\lt,  \quad & [E_\beta]_{ij} &= \beta(\zeta_i, \zeta_j)_\lt,
\end{aligned}
\end{equation}
for each basis function ${v}_i \in U_{0,\hp}$ and $\zeta_i \in \Psi_\hp$ and user-chosen parameter $\beta \geq 0$. $E_\beta$ is augmented to the bottom right block in order reduce the condition number of $G$ but is not required for well-posedness. In fact, in most examples in \Cref{sec:examples}, we fix $\beta=0$ which corresponds to $E_\beta \equiv 0$. Note that $D_\psi$ is the only matrix with dependence on the current Newton iterate.

\subsection{Gradient-type constraint}
\label{sec:newton:gradient}

In order to discretize the gradient-type LVPP subproblem \cref{eq:pG:gradient} we choose a similar discretization. Now we discretize $u$ with a (partial) degree $\vect{p}$  continuous $p$-FEM basis and $\psi$ with the (partial) degree $\vect{p}-1$ discontinuous Legendre basis, for $\vect{p}\geq 1$. This choice of pair is motivated by \cite[Ex.~6]{dokken2024}. The numerical analysis is less developed for this problem and other inf-sup stable pairs may exist.

After a Newton linearization, this leads to the same symmetric saddle point structure  as \cref{eq:pG-matrix} that we choose not to relabel. The $\alpha$-scaled stiffness matrix $A_\alpha$ is identical to the one in the previous subsection. Recall that $\Phi_\hp$ denotes the discontinuous basis of Legendre spectral Galerkin basis defined in \Cref{sec:fem}. Then, for each basis function ${v}_i \in U_{0,\hp}$, $\zeta_i \in \Psi_\hp^d$, and $\eta_i \in \Phi_\hp^d$, $d \in \{1,2\}$, the (not relabeled) matrices $B$, $D_\psi$, and $E_\beta$ are defined as:
\begin{equation}
\begin{aligned}
B_{ij} =  (\nabla v_i, \zeta_j)_\lt, \;\;  [D_\psi]_{ij} = \langle F'_{\zeta_i}(\psi_\hp), \zeta_j \rangle,  \;\; [E_\beta]_{ij} = \beta (\nabla_h \eta_i, \nabla_h \eta_j)_\lt,
\end{aligned}
\end{equation}
where $\nabla_h$ denotes the broken (cellwise) gradient,  $F'_{\zeta_i}(\psi_\hp)$ represents the Fr\'echet derivative of $F_{\zeta_i}(\psi_\hp) \coloneqq (\zeta_i, \varphi_\hp \psi_\hp (1+|\psi_\hp|^2)^{-1/2})_{L^2(\Omega)}$, and $\varphi_\hp$ denotes the FEM approximation of $\varphi$. $E_\beta$ plays the same role as in \Cref{sec:newton:obstacle} and, once again, $D_\psi$ is the only matrix with dependence on the current Newton iterate. The effectiveness of the hpG algorithm hinges on the ability to solve \cref{eq:pG-matrix} efficiently and in a manner that is robust to choices of the mesh size $h$, the truncation degree $p$, $\alpha$, and $\psi_\hp$.

\subsection{The matrices $A_\alpha$, $B$, $D_\psi$, and $E_\beta$}
The matrices $A$, $B$, and $E_\beta$ are independent of any of the parameters of the pG algorithm and sparse with $O(p^d/h^d)$ entries that are nonzero. In particular $E_\beta$ is block-diagonal after a permutation of the rows and columns. Hence they can be efficiently assembled, stored and applied to a vector in $O(p^d/h^d)$ flops. Moreover, when $d=1$, $A$ admits a reverse Cholesky factorization\footnote{A \emph{reverse} Cholesky factorization $A = L^\top L$ initializes the factorization from the bottom right corner of a symmetric positive-definite matrix $A$ rather than the top left.} $A=   L^\top L$ such that the factor $L$ is also sparse and has $O(p/h)$ nonzero entries \cite[Th.~4.2]{knook2024quasi}. Thus the inverse of $A_\alpha$ may be computed and applied with $O(p/h)$ flops. In two dimensions, we observe that the Cholesky factorization of $A$ has only a moderate fill-in and can be computed quickly. In this paper, we always performed this factorization once at the start of a solve. However, we note that an optimal complexity iterative solver also exists for $A$ in 2D \cite{knook2024quasi,papadopoulos2024b} via the ADI algorithm \cite{Fortunato2020}.

The matrix $D_\psi$ is not as sparse as $A_\alpha$ and $B$. However, since all the basis functions in $\Psi_\hp$ are supported on a maximum of one cell, the columns and rows of $D_\psi$ may be permuted to a block diagonal structure although the blocks themselves may be dense, leading to $O((p-j)^{2d}/h^d)$ nonzero entries where $j=2$ for the obstacle problem and $j=1$ in the gradient-type constraint case. We note that we can apply the action of $D_\psi$ with quasi-optimal $O((p-j)^d \log^d (p-j) / h^d )$ flops. Since we are able to apply $D_\psi$ in an efficient manner and have a fast inverse for the top left block $A_\alpha$, this motivates a block preconditioning approach coupled with a matrix-free Krylov method. 

\begin{figure}[h!]
\centering
\subfloat[$A_\alpha$.]{\includegraphics[width =0.19\textwidth]{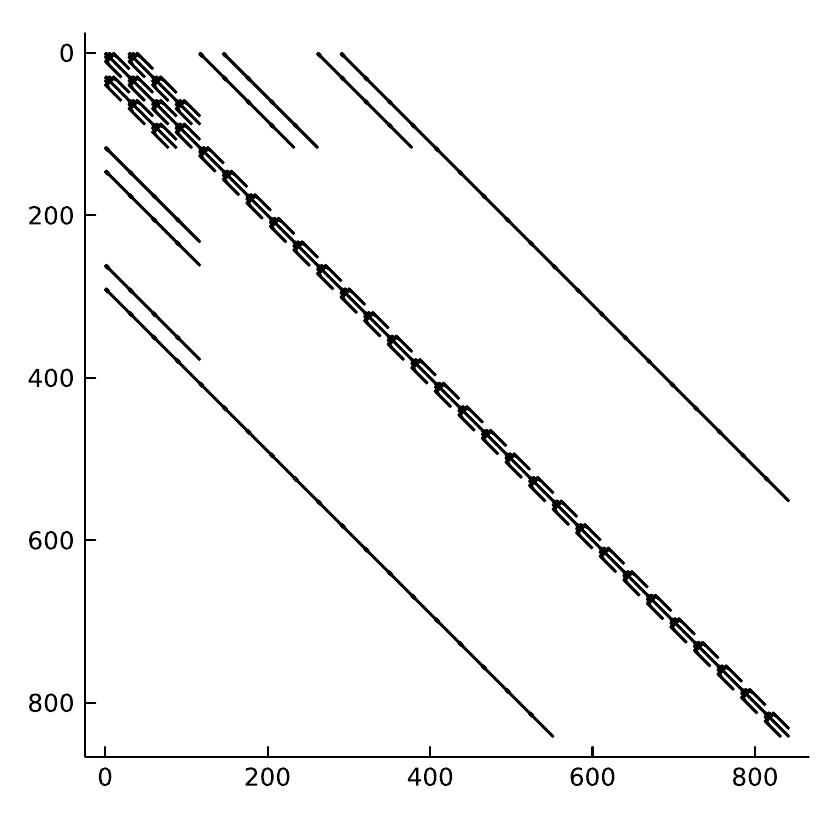}}
\subfloat[$L^{\text{chol}}_\alpha$.]{\includegraphics[width =0.19\textwidth]{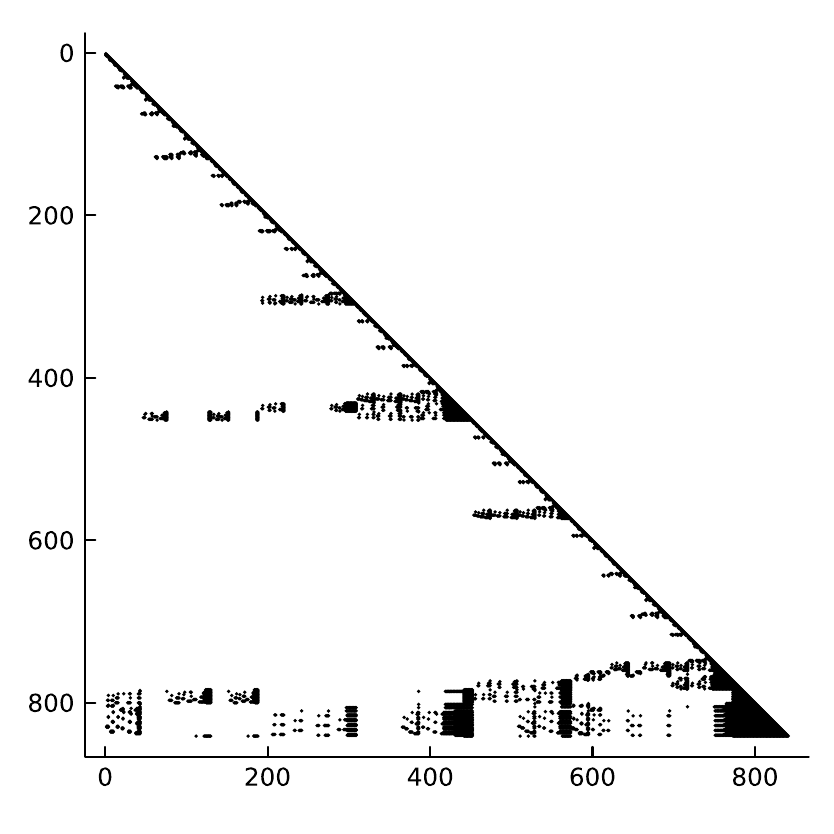}}
\subfloat[$B$.]{\includegraphics[width =0.19\textwidth]{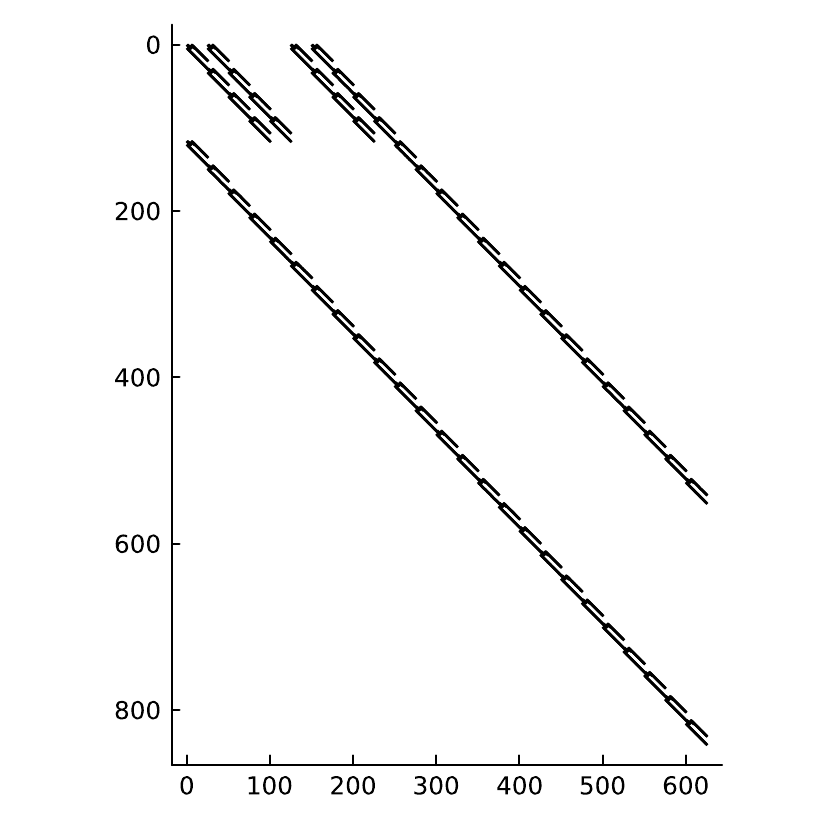}}
\subfloat[$D_\psi$.]{\includegraphics[width =0.19\textwidth]{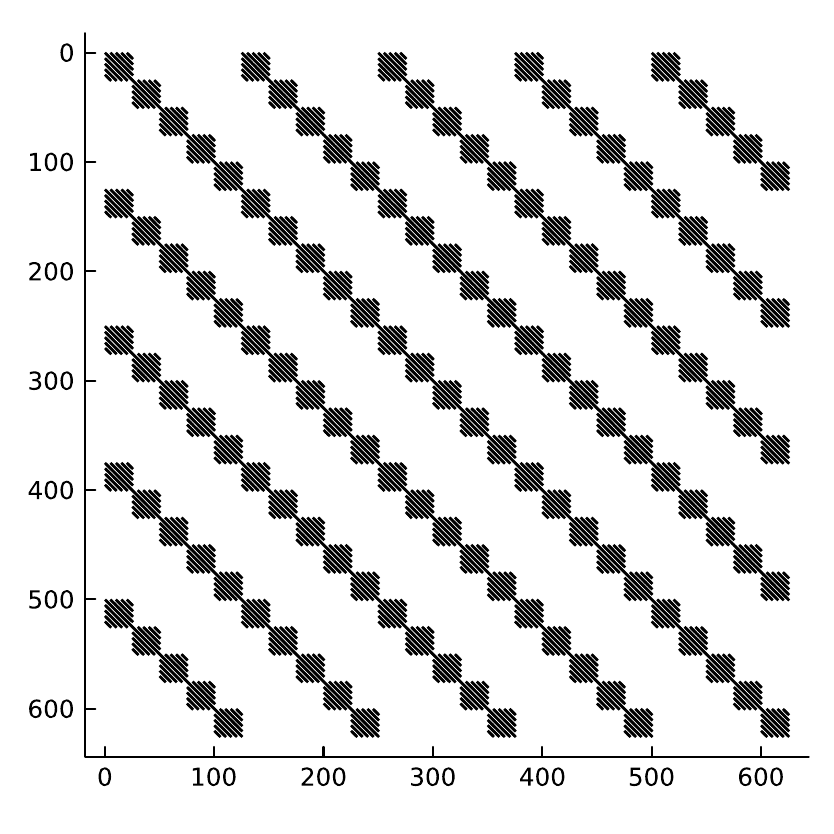}}
\subfloat[Permuted $D_\psi$.]{\includegraphics[width =0.19\textwidth]{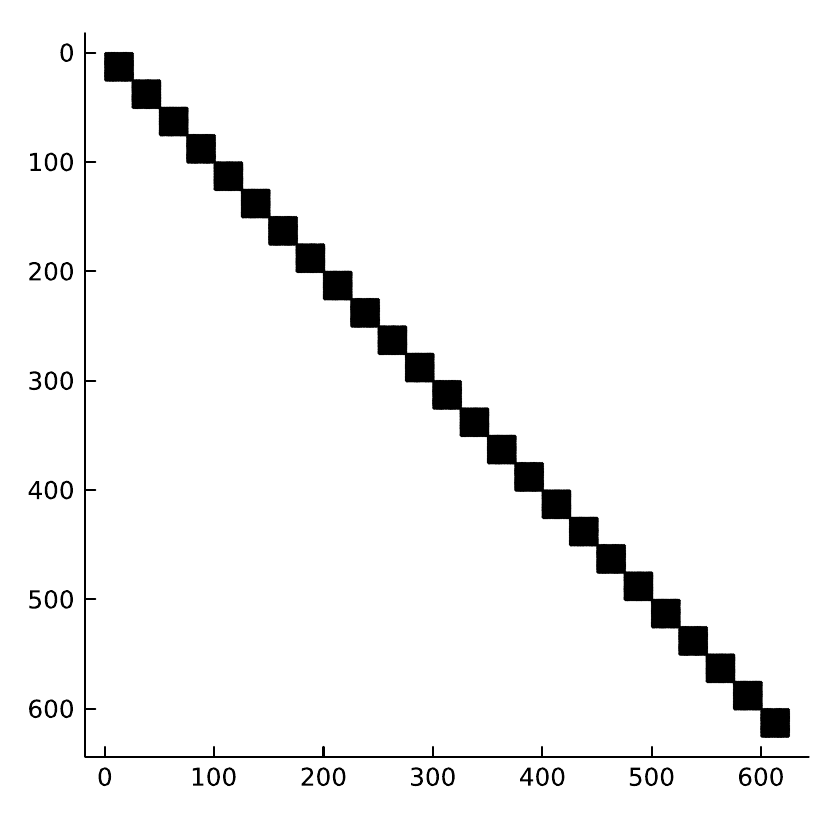}}
\caption{The spy plots of the scaled stiffness matrix $A_\alpha$, its lower Cholesky factor $L^{\text{chol}}_\alpha$ or reverse lower Cholesky factor $L^{\text{rchol}}_\alpha$, the Gram matrix $B$, the exponential block $D_\psi$, and $D_\psi$ with permuted rows and columns to reveal the hidden block diagonal structure for the Newton linearization of the pG obstacle subproblem \cref{eq:pG} when using a two-dimensional tensor-product $p$-FEM discretization with 25 cells and partial degree $p=5$ on each cell.}
\label{fig:spy-plots}
\end{figure}

In \cref{fig:spy-plots} we provide spy plots of the individual blocks in $G$ in a two-dimensional discretization for the pG obstacle subproblem \cref{eq:pG}. The gradient-type case results in similar sparsity patterns. The discretization involves 25 cells with partial degree $p=5$ on each cell. Note the sparsity in $A_\alpha$, its (reverse) Cholesky factor, and $B$. Moreover, $D_\psi$ has a block diagonal structure (after rearranging the rows and columns) and hence computing its inverse can be performed cellwise and in parallel.

\subsection{Block preconditioning}
\label{sec:preconditioning}

Given that the submatrices within $G$ in \cref{eq:pG-matrix} are relatively sparse, one can assemble the full matrix and perform a sparse LU factorization--a strategy that showed unexpectedly good performance in our examples. However, by employing a straightforward block preconditioning strategy, we can improve the conditioning of the solver and achieve faster solve times, particularly as $h \to 0$ and $p \to \infty$. Block preconditioning reduces the problem to solving a smaller system involving a dense Schur complement matrix. Despite its density, this matrix allows for an efficient matrix-vector product, making it well-suited for solution via Krylov subspace methods. As usual, maintaining an acceptable iteration count depends on the availability of a suitable preconditioner. In this subsection we develop a (up to permutation) block-diagonal  preconditioner, whose effectiveness is expected to deteriorate at a worse-case polylogarithmic rate  with respect to $p$ and $h$. %

Via a Schur complement factorization of the matrix $G$ in \cref{eq:pG-matrix}, one finds the inverse of $G$ admits the triple-matrix-product decomposition \cite[Sec.~5]{benzi2005}:
\begin{align}
G^{-1} = P_R^{-1} P_D^{-1} P_L^{-1} =
\begin{pmatrix}
I & -A_\alpha^{-1} B \\0 & I
\end{pmatrix}
\begin{pmatrix}
A_\alpha^{-1} & 0 \\ 0 & S^{-1}
\end{pmatrix}
\begin{pmatrix}
I & 0 \\ -B^\top A_\alpha^{-1} & I
\end{pmatrix},
\label{eq:schur4}
\end{align}
where the Schur complement matrix $S$ is $S \coloneqq -(D_\psi +E_\beta+ B^\top A_\alpha^{-1} B)$. If one has access to (approximations of) $A_\alpha^{-1}$ and $S^{-1}$, then one may invert $G$ via $P_L$, $P_D$, and $P_R$. 

\begin{figure}[ht]
\centering
 \begin{tikzpicture}[%
every node/.style={draw=black, thick, anchor=west},
grow via three points={one child at (-0.5,-0.7) and
    two children at (-0.5,-0.7) and (-0.5,-1.4)},
edge from parent path={(\tikzparentnode.190) |- (\tikzchildnode.west)}]
\node {Outer Krylov solver for $G$ (FGMRES)}
child { node {Block preconditioner ($P=P_D, P_L P_D$ or $P_F$)}
    child { node {Invert $A_\alpha$ (Cholesky or preconditioned CG)}}
        child { node{Inner Krylov solver for $S$ (GMRES)}
        	 child { node{Block-diagonal preconditioner $\hat{S}$}} 
    }
};
\end{tikzpicture}
\caption{The solver diagram for inverting $G$ in \cref{eq:pG-matrix}.}\label{fig:solver}
\end{figure}

In \cref{fig:solver} we outline the building blocks for an iterative solver for $G$. First, one leverages an outer FGMRES solver \cite{saad1993}.\footnote{FGMRES, rather than GMRES, must be used if an inner GMRES solver is used to apply $S^{-1}$.} For this to converge in a reasonable number of iterations we require a suitable preconditioner. The Schur complement reduction in \cref{eq:schur4} suggests three natural block preconditioners. One may choose the full factorization $P = P_F \coloneqq P_L P_D P_R$. With exact inverses of $A_\alpha$ and $S$, this is the perfect preconditioner since $P_R^{-1} P_D^{-1} P_L^{-1} G = I$. Hence, the outer FGMRES solver would converge in one iteration. Naively, one might think utilizing $P_F$ requires three applications of $A_\alpha^{-1}$ and one of $S^{-1}$, however, applying $P_F^{-1}$ is equivalent to the following sequential solve procedure:
\begin{align}
\vectt{y} = \vectt{b}_\psi - B^\top A_\alpha^{-1} \vectt{b}_u, \quad \vectt{\delta}_\psi= S^{-1} \vectt{y}, \quad
\vectt{\delta}_u = A_\alpha^{-1} (\vectt{b}_u - B \vectt{\delta}_\psi). \label{eq:schur2}
\end{align}
Hence, one only requires two applications of $A_\alpha^{-1}$ and one of $S^{-1}$. The two other common choices of block preconditioner are $P = P_D$ and $P = P_L P_D$ \cite[Sec.~4]{pearson2017}. Both of these choices only require one application of $A_\alpha^{-1}$ and one of $S^{-1}$.

As is typical with the Schur complement, $S$ is dense due to the $B^\top A_\alpha^{-1} B$ term. Hence, it is computationally catastrophic to assemble.  Irrespective of the choice for $P$, we must construct an approximation for $S$ that admits a cheap assembly and solve. In the subsequent sections, we shall build a block-diagonal approximation $\hat{S}$. If the action of $A_\alpha^{-1}$ may computed efficiently, then the action of $S$ can be applied efficiently (each matrix in the triple matrix product $B^\top A_\alpha^{-1} B$ is applied sequentially). This motivates utilizing a $\hat{S}$-preconditioned inner GMRES solver to apply $S^{-1}$ to any user-chosen tolerance. If the tolerance is sufficiently small, then the outer FGMRES solver will converge in one iteration if $P=P_F$. %

The choice of $P$ is inherently determined by the cost of computing $A_\alpha^{-1}$. If it is cheap, one should choose $P = P_F$, since out of the three choices, this will reduce the number of outer Krylov iterations and hence the number of required approximate actions of $S^{-1}$. One option is to perform and cache a single Cholesky factorization of $A$ at the very beginning of the outer nonlinear hpG solve.  Thanks to the high-order discretizations, the matrix $A$ remains much smaller whilst achieving the same error when compared to a $p=1$ discretization. Whereas the stiffness matrix $A$ of a $p=1$ discretization may quickly become sufficiently large so that direct factorizations are infeasible, this bottleneck occurs much later in our $hp$-FEM framework. We observe this phenomenon in all the one- and two-dimensional examples of \Cref{sec:examples} where the Cholesky factorization is attained in fractions of a second and its application in tens of times faster than that. Thus choosing $P=P_D$ or $P=P_L P_D$ in these situations leads to a slower solve, since we must apply the slower action $S^{-1}$ more times. Hence, in all two-dimensional examples, we shall solely consider the solver strategy of \cref{fig:solver2}.
\begin{figure}[ht]
\centering
 \begin{tikzpicture}[%
every node/.style={draw=black, thick, anchor=west},
grow via three points={one child at (-0.5,-0.7) and
    two children at (-0.5,-0.7) and (-0.5,-1.4)},
edge from parent path={(\tikzparentnode.200) |- (\tikzchildnode.west)}]
\node {Invert $G$}
child { node {Schur complement factorization}
    child { node {Invert $A_\alpha$ (Cholesky)}}
        child { node{Krylov solver for $S$ (GMRES)}
        	 child { node{Block-diagonal preconditioner $\hat{S}$}} 
    }
};
\end{tikzpicture}
\caption{The solver diagram for inverting $G$ when $d=2$.}\label{fig:solver2}
\end{figure}
However, in three dimensions, the size of $A_\alpha$ grows sufficiently quickly that a Cholesky factorization becomes infeasible sooner. Without a factorization, one must approximate the action of $A_\alpha^{-1}$ at each Krylov iteration and, therefore, may expect savings if fewer applications are required. We explore this route in \Cref{sec:examples:3d}, where we approximate the action of $A_\alpha^{-1}$ via algebraic multigrid preconditioned CG and compare the solve times whilst using $P=P_F, P_LP_D$ or $P=P_D$ as well as changing the tolerance of the inner Krylov solver to limit the number of applications of $S$ (and subsequently $A_\alpha^{-1}$). 

In the next two subsections we design a block-diagonal approximation $\hat{S}$ to $S$.

\subsection{Obstacle problem: a Schur complement preconditioner}
\label{sec:preconditioning:obstacle}
In this subsection we focus on preconditioning the Schur complement as it arises in the context of the pG obstacle subproblem \cref{eq:pG}. For ease of implementation, we will restrict our search for a preconditioner for $S$ in \cref{eq:schur2} that can be assembled and factorized. Therefore, we do not consider preconditioners based on $p$-multigrid in this work \cite{Babuska1991,beuchler2006b,korneev1999,Schwab1998}. 

Our goal is to approximate $B$ and $A_\alpha$ with the alternative matrices $\hat{B}$ and $\hat{A}_\alpha$ such that the new matrix-triple-product $\hat{B}^\top \hat{A}^{-1}_\alpha \hat{B}$ induces a (up to permutation) block-diagonal matrix. Common choices involve removing or modifying the hat functions in $U_{0,\hp}$ \cite[Sec.~3]{Babuska1991}. We opt for removing the hat functions and rescaling the remaining shape functions, i.e.~we use the Legendre spectral Galerkin basis $\Phi_\hp$ defined in \Cref{sec:fem} which as we found it exhibited the lowest iteration counts. We choose the following sparse approximation for $S$
\begin{align}
\hat{S} = -D_\psi- E_\beta -  \hat{B}^\top \hat{A}_\alpha^{-1} \hat{B}
\label{eq:schur-preconditioner}
\end{align}
where, for each basis function $\eta_i \in \Phi_\hp$ and $\zeta_i \in \Psi_\hp$,
\begin{align}
[\hat{A}_\alpha]_{ij} \coloneqq \alpha (\nabla_h \eta_i, \nabla_h \eta_j)_\lt \;\; \text{and} \;\; \hat{B}_{ij} \coloneqq (\eta_i, \zeta_j)_\lt.
\end{align}
$\hat{A}_\alpha$ is diagonal when $d=1$ and block-diagonal  (after a permutation of the rows and columns) with 4N blocks (where $N$ is the number of 2D cells in the mesh) when $d=2$. Thus $\hat{A}$ admits a cheap sparse Cholesky factorization. $\hat{S}$ is also block-diagonal although the blocks themselves are dense and of size $O(p^d \times p^d)$, $d \in \{1,2,3\}$.

\begin{figure}[h!]
\centering
\subfloat[16 quad.~cells]{\includegraphics[width =0.24\textwidth]{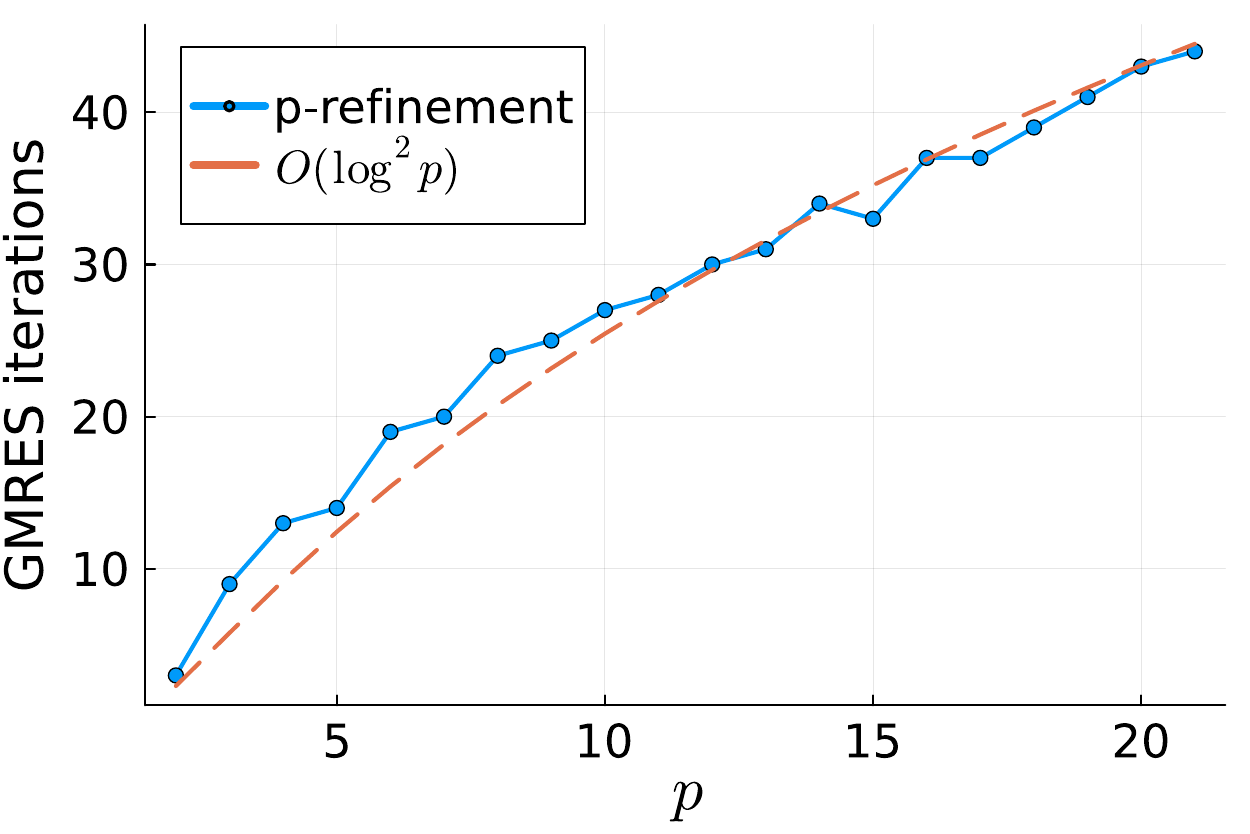}}
\subfloat[16 quad.~cells]{\includegraphics[width =0.24\textwidth]{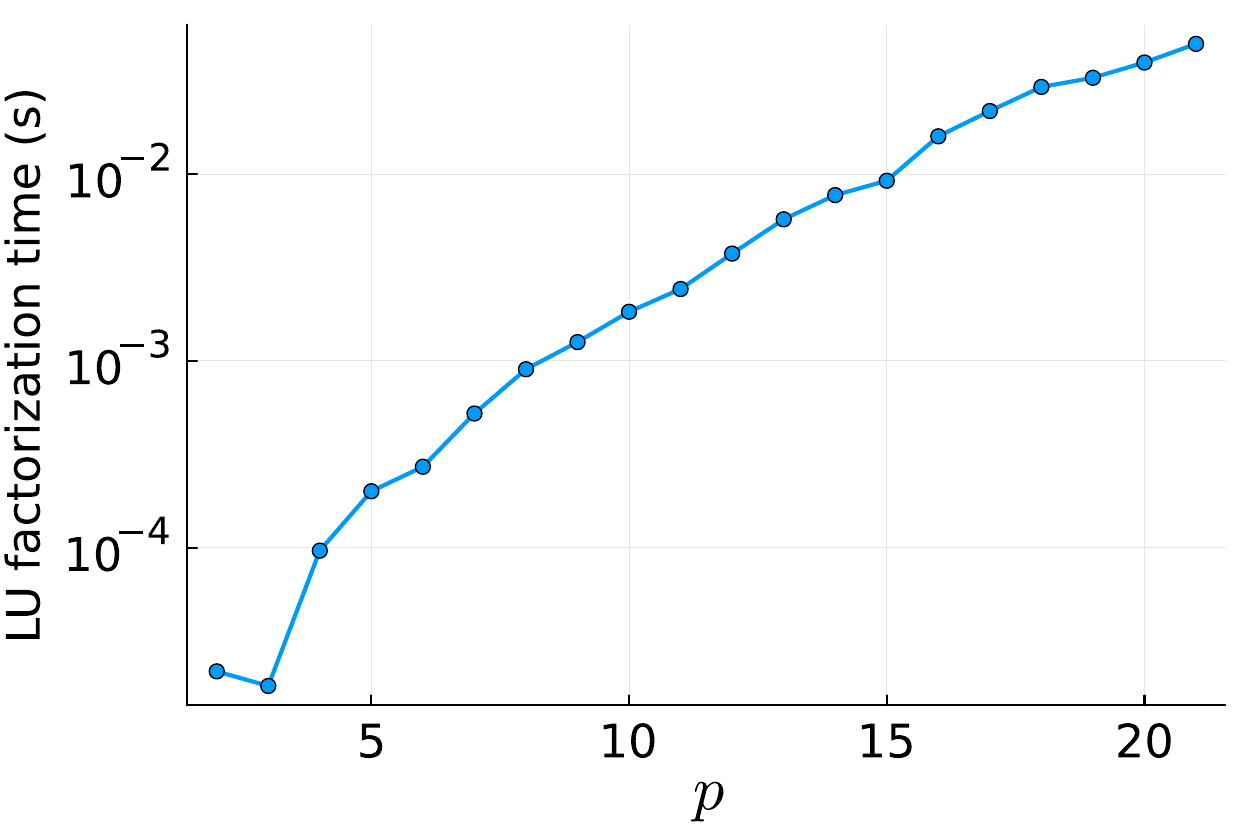}}
\subfloat[$p=5$]{\includegraphics[width =0.24\textwidth]{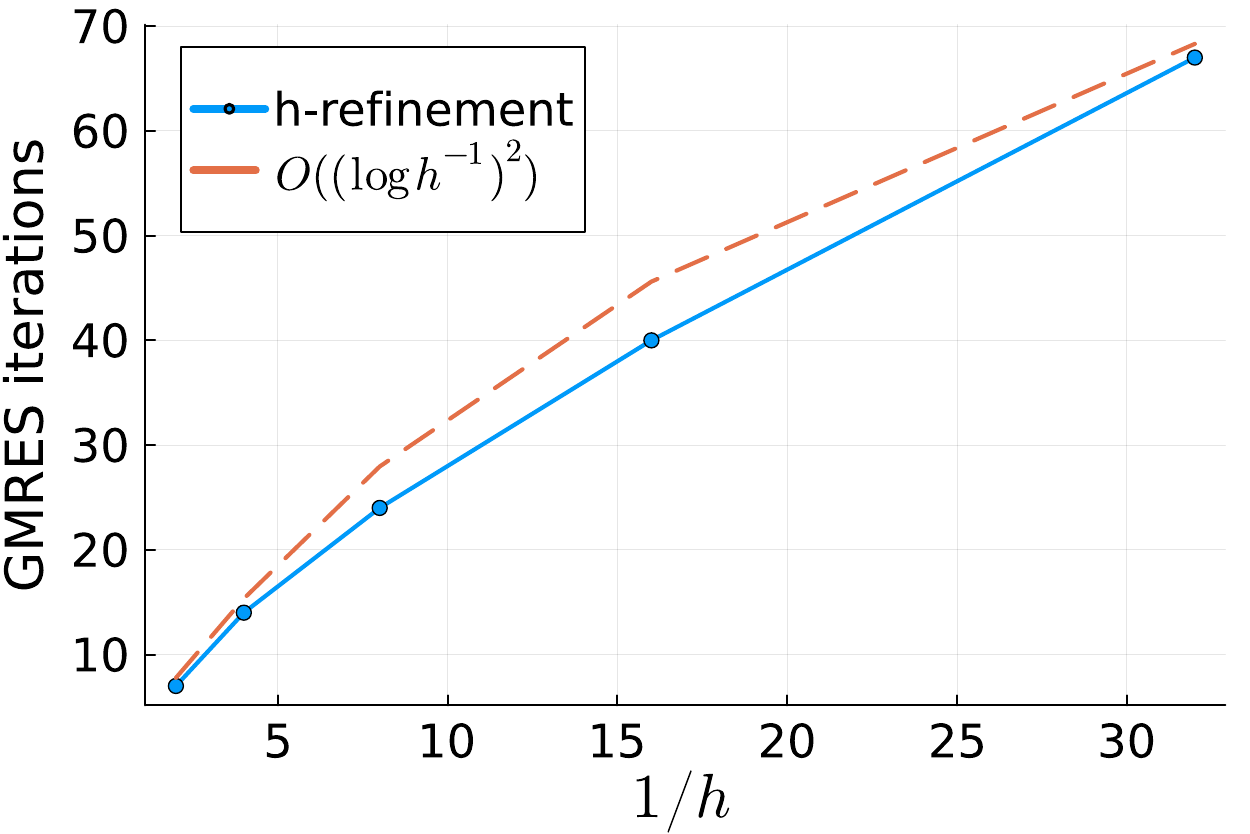}}
\subfloat[$p=5$]{\includegraphics[width =0.24\textwidth]{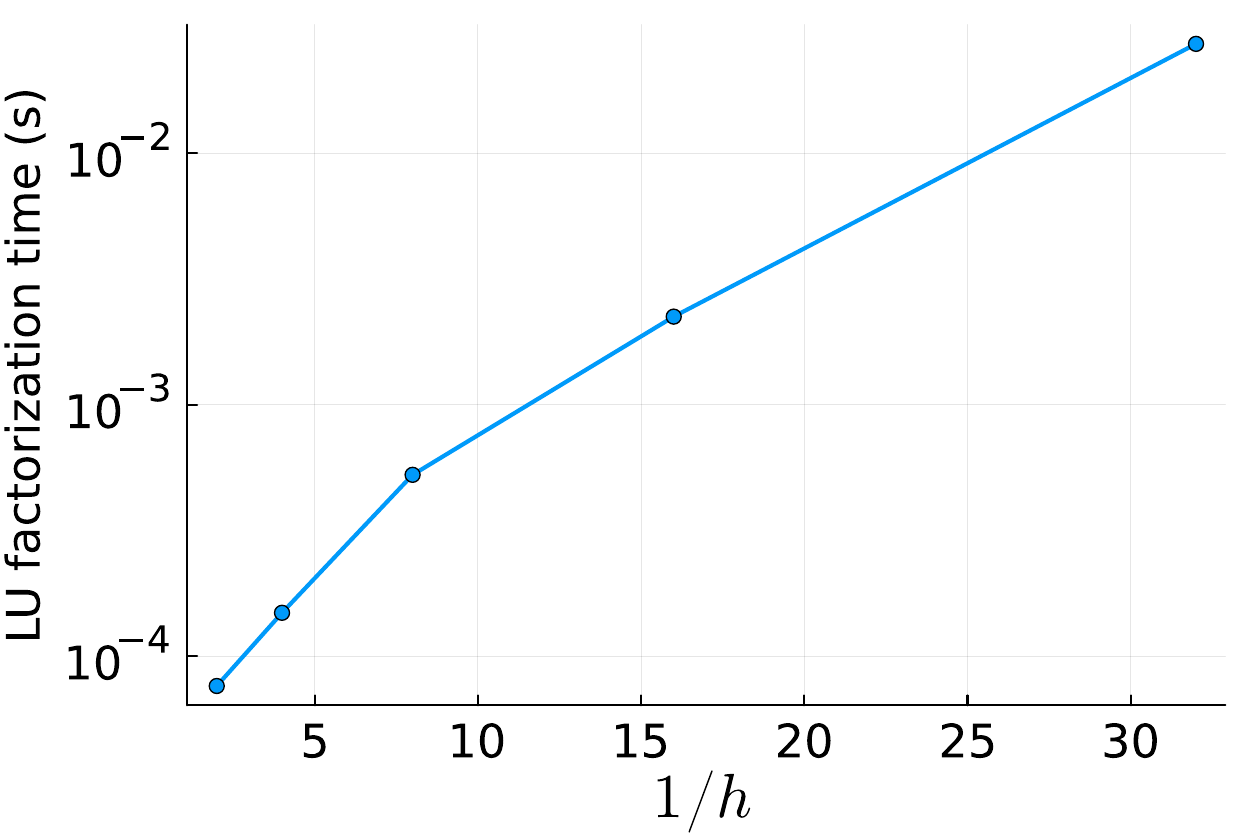}}
\caption{(Obstacle problem). We consider $d= 2$ and fix $\alpha =1$, $\beta=0$, $D_\psi \equiv 0$ and consider the dense Schur complement $S$ as defined \cref{eq:schur2}, in the context of the obstacle subproblem \cref{eq:pG}, and its sparse preconditioner $\hat{S}$ as defined in \cref{eq:schur-preconditioner}. We measure the growth in the number of $\hat{S}$-left-preconditioned GMRES iterations to solve $S \vectt{x} = \vect{1}$ to a relative error of $10^{-6}$ (i.e.~$\|\vectt{r}_k\|_{\ell^2}/\|\vectt{r}_0\|_{\ell^2} \leq 10^{-6}$) with respect to partial degree $p$ in (a) and $1/h$ in (c). We observe polylogarithmic growth. We  measure the LU factorization time for $\hat{S}$ with respect to partial degree $p$ in (b) and $1/h$ in (d).}
\label{fig:preconditioner}
\end{figure}

In \cref{fig:preconditioner} we consider \cref{eq:pG-matrix} as constructed for a 2D discretization and measure the growth of $\hat{S}$-left-preconditioned GMRES iterations to solve $S \vectt{x} = \vect{1}$ to a relative error of $10^{-6}$ with respect to partial degree $p$ and $1/h$. We fix $\alpha=1$, $\beta=0$, and $D_\psi \equiv 0$ which constitutes the worst case scenario, i.e.~where $S$ is the most singular. We observe polylogarithmic growth in the number of GMRES iterations. We also provide timings for the sparse LU factorization of $\hat{S}$.

\subsection{Gradient-type constraints: a Schur complement preconditioner}
\label{sec:preconditioning:gradient}
\begin{figure}[h!]
\centering
\subfloat[16 quad.~cells]{\includegraphics[width =0.24\textwidth]{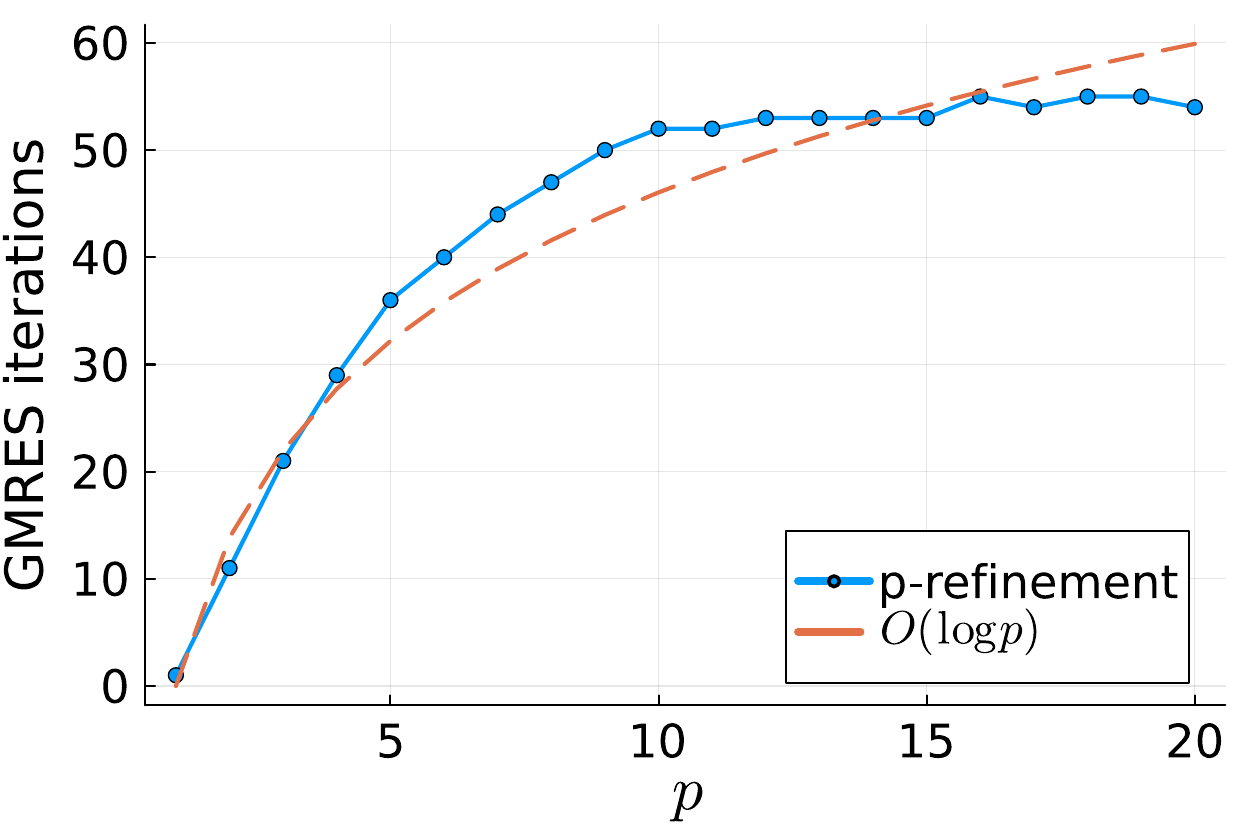}}
\subfloat[16 quad.~cells]{\includegraphics[width =0.24\textwidth]{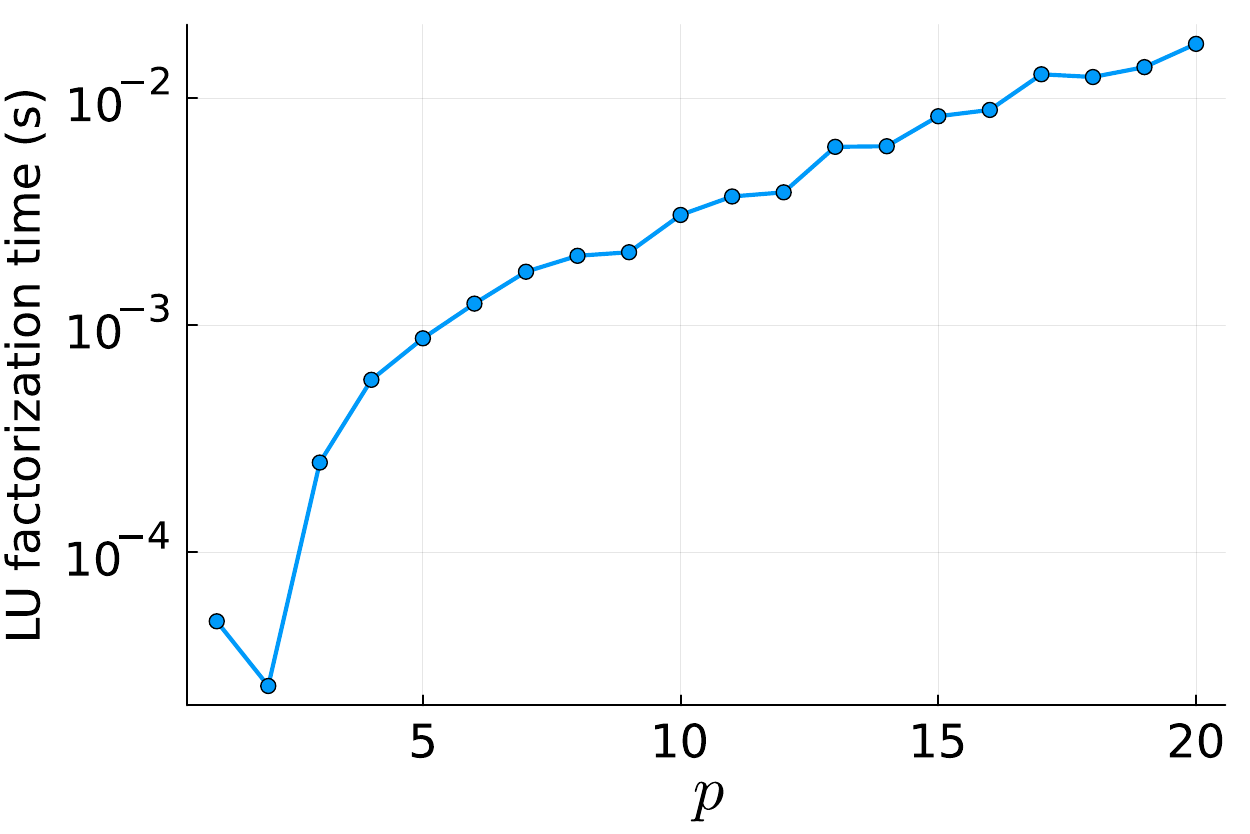}}
\subfloat[$p=4$]{\includegraphics[width =0.24\textwidth]{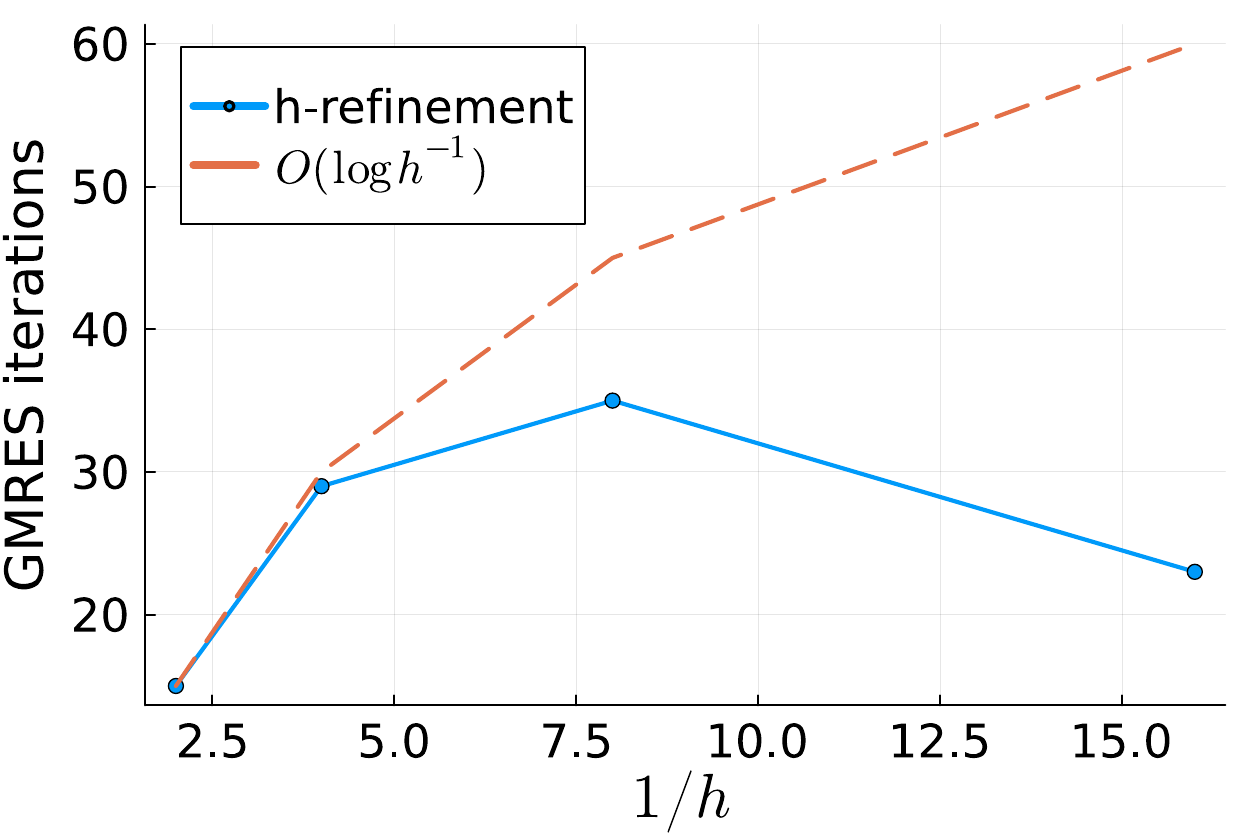}}
\subfloat[$p=4$]{\includegraphics[width =0.24\textwidth]{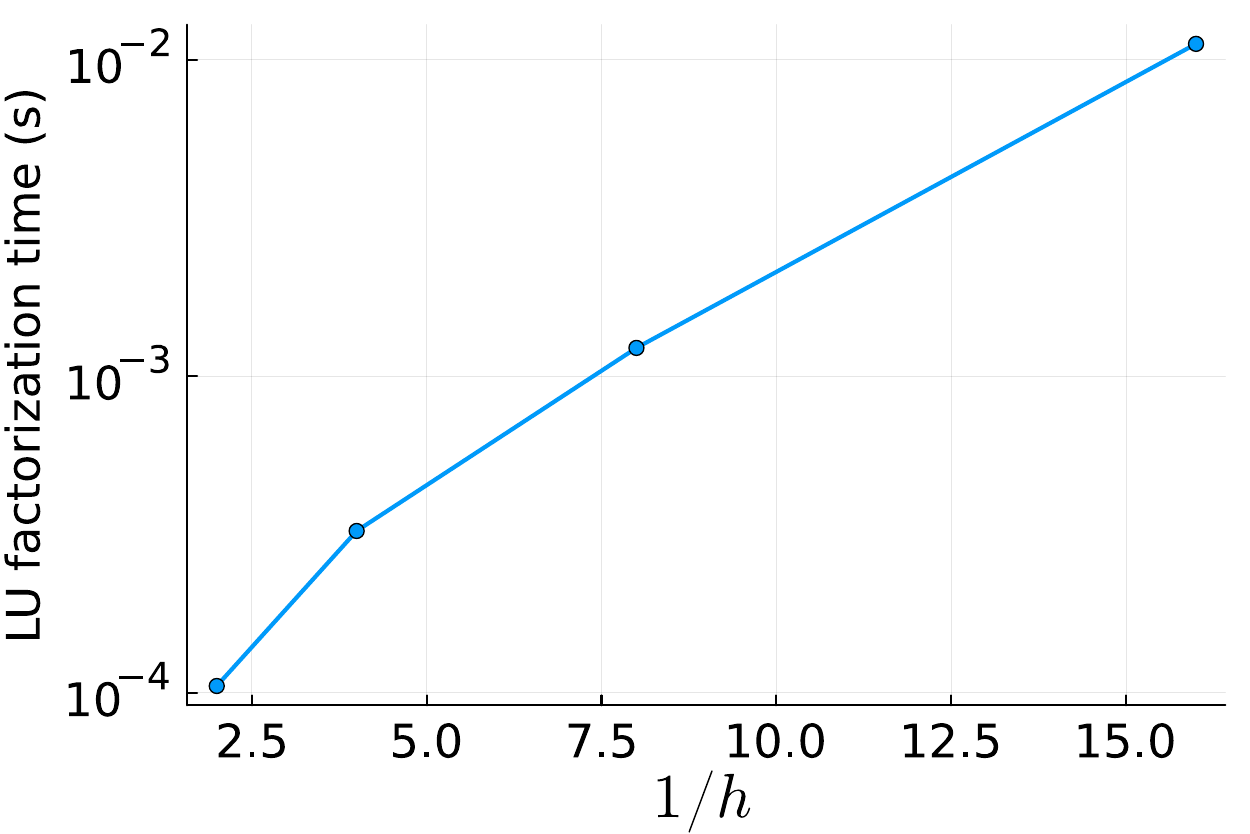}}
\caption{(Gradient-type constraint). We consider $d=2$ and fix $\alpha =1$, $\beta=10^{-5}$, $D_\psi \equiv 0$ and consider the dense Schur complement $S$ as defined \cref{eq:schur2}, in the context of the pG gradient-type subproblem \cref{eq:pG:gradient}, and its sparse preconditioner $\hat{S}$ as defined in \cref{eq:schur-preconditioner:gradient}. We measure the growth in the number of $\hat{S}$-left-preconditioned GMRES iterations to solve $S \vectt{x} = \vect{1}$ to a relative error of $10^{-6}$ with respect to partial degree $p$ in (a) and $1/h$ in (c). We observe slower than logarithmic growth in both cases. We also measure the sparse LU factorization time for $\hat{S}$ with respect to partial degree $p$ in (b) and $1/h$ in (d).}
\label{fig:preconditioner:gradient}
\end{figure}

We choose a simpler preconditioner in the case of the pG gradient-type subproblem \cref{eq:pG:gradient}. More specifically, we drop the triple-matrix-product in $S$ and choose:
\begin{align}
\hat{S} = -D_\psi - E_\beta.
\label{eq:schur-preconditioner:gradient}
\end{align}
This provided robust preconditioning as evidenced in \cref{fig:preconditioner:gradient}. We perform the same experiment as for the obstacle-type case and examine the number of left-preconditioned GMRES iterations to solve $S\vectt{x} = \vectt{1}$ to a relative tolerance of $10^{-6}$ with respect to increasing $p$ and $1/h$ in 2D. We fix $\alpha=1$, $D_\psi \equiv 0$ but select $\beta = 10^{-5}$ as otherwise the Schur complement $S$ is nearly singular. 

\section{Notes for fast implementation}
\label{sec:implementation}

\indent\indent \textbf{Fast transforms.} An essential aspect for fast solve times when handling high-degree polynomials is the ability to implement the analysis (expansion in the basis) and synthesis (evaluation of the basis on a grid) routines in quasi-optimal complexity. Such transforms exist for the $L^2$-conforming basis $\Psi_\hp$ consisting of scaled-and-shifted Legendre polynomials \cite{alpert1991, keiner2011, townsend2018, SOActa}. In our implementation we opt for the approach described in \cite[Sec.~4.6.4]{SOActa} and implemented in \cite{fasttransforms.jl} where the coefficients of the Legendre expansion are converted to and from the coefficients of the equivalent Chebyshev expansion. The Chebyshev expansion then enjoys the use of the Discrete Cosine Transform (DCT) for quasi-optimal analysis and synthesis.

\indent \textbf{Alleviating ill-conditioning.} The nonlinear convergence of the Newton solver may degrade as $h \to 0$ and $p \to \infty$ due to numerical ill-conditioning of the linear system in \cref{eq:pG-matrix}. We find that there are three remedies to such a situation:
\itemsep=-2pt
(1) Utilization of the block preconditioning strategy in \Cref{sec:preconditioning}.
(2) Adding a small modification to the Jacobian, which we denoted by $E_\beta$ in \cref{eq:pG-matrix}. The local modification also makes the GMRES solver more robust.
(3) Choosing small values for the $\alpha$-sequence. Since the solver converges without requiring $\alpha \to \infty$, one may use small values of $\alpha$, for instance terminating whilst $\alpha < 1$, and still observe numerical convergence.

\indent \textbf{Quadrature.} Quadrature is required to assemble or apply the action of the matrix $D_\psi$. For instance, in the obstacle problem, given $\psi_\hp$ we must compute (I) $(\zeta_i, \E^{-\psi_\hp})_{L^2(\Omega)}$ and (II) $(\zeta_i, \E^{-\psi_\hp} \zeta_j)_{L^2(\Omega)}$ for each basis function $\zeta_i, \zeta_j \in \Psi_\hp$. In our implementation we opted for a fast but non-standard quadrature. We consider the approximations $\E^{-\psi_\hp} \approx \sum_{k} c_k \zeta_k$ for (I) and $\E^{-\psi_\hp} \zeta_j \approx \sum_{k} \tilde{c}_k \zeta_k$ for (II). Once the expansions are performed, the integrals (I) and (II) are computed efficiently due to the orthogonality and support properties of the basis functions in $\Psi_\hp$. The integrals are only exact if $\E^{-\psi_\hp}$ is a piecewise constant. Nevertheless, experimentally the quadrature rule provided excellent results with the hpG solver exhibiting robustness against polynomial aliasing \cite[Ch.~2.4.1.2]{Karniadakis2005}. This quadrature scheme does not preserve symmetry in $D_\psi$ but when tested against a symmetry-conserving Clenshaw--Curtis quadrature scheme \cite{Clenshaw1960, trefethen2008} we found that the non-standard route allowed for quicker assembly of $D_\psi$ and actually reduced the number of nonlinear and GMRES iterations required for the obstacle problem solve. We use the same technique for the action and assembly of $D_\psi$ in the gradient-type case.

\indent \textbf{Inexact solves.} Speedups are available if one does not solve the Newton linear systems exactly but rather terminates the linear system  iterative solver once the error falls below a user-chosen tolerance. Theoretical guarantees are an open problem.

\section{Examples}
\label{sec:examples}

\textbf{Data availability.} 
The implementation is contained in the package HierarchicalProximalGalerkin.jl \cite{hier-obstacles.jl} and is written in the open-source language Julia \cite{Bezanson2017}. The version of  HierarchicalProximalGalerkin.jl run in our experiments is archived on Zenodo \cite{hpg-zenodo}. We utilize a number of available Julia packages \cite{linesearches.jl, Badia2020, extendablesparse.jl, iterativesolvers.jl, linearmaps.jl, matrixfactorizations.jl, fasttransforms.jl} and, in particular, the
PiecewiseOrthogonalPolynomials.jl package \cite{piecewisepolys.jl} for its implementation of the hierarchical $p$-FEM basis. The experiments in \Cref{sec:examples:1d,sec:examples:2d:1,sec:examples:gradient,sec:examples:3d} were run on a deskstop with 16GB of RAM and 8 CPUs Intel(R) Core(TM) i7-10700 CPU @ 2.90GHz. The penultimate example in \Cref{sec:examples:qvi} was run on a machine with 768GB of RAM and 72 CPUs HPE Synergy 660 Gen10
Xeon @ 3.1GHz.

\subsection{Oscillatory data}
\label{sec:examples:1d}
\begin{figure}[h!]
\centering
\includegraphics[width =0.45 \textwidth]{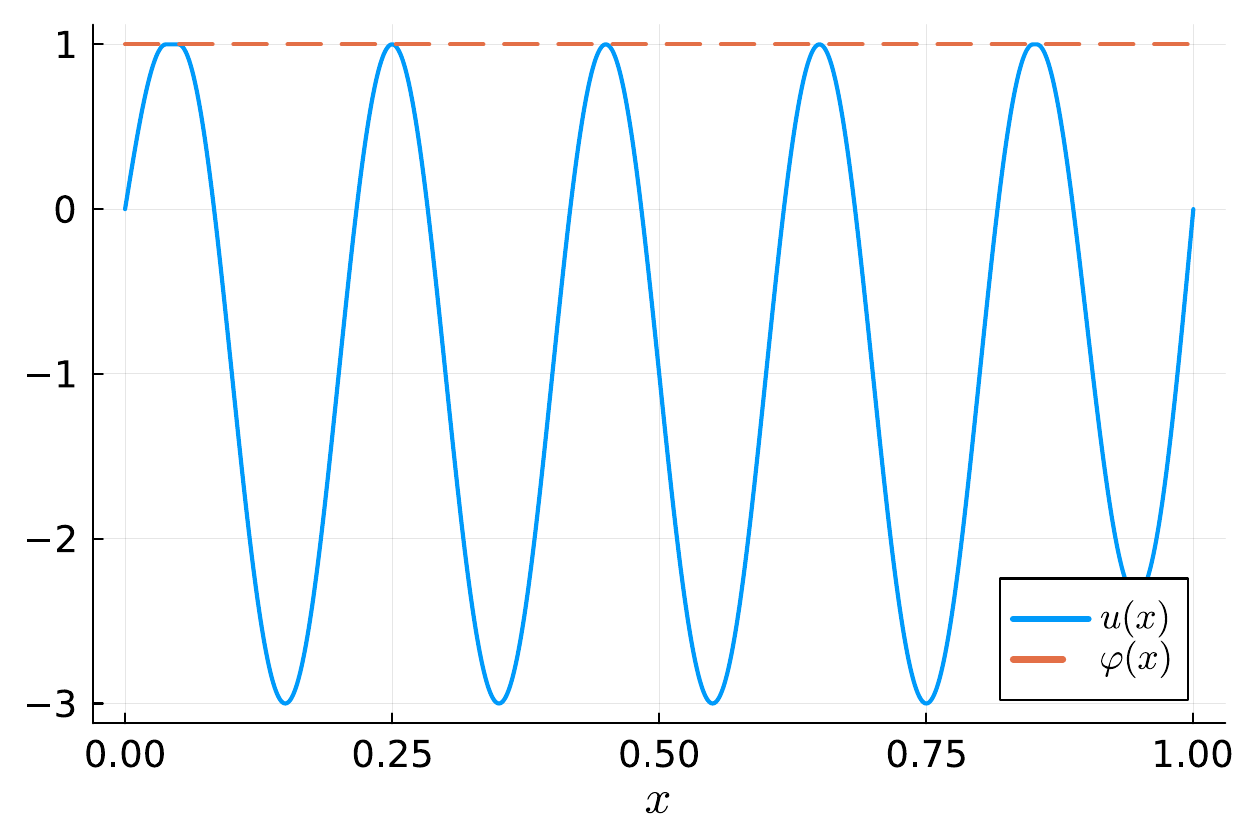}
\includegraphics[width =0.45 \textwidth]{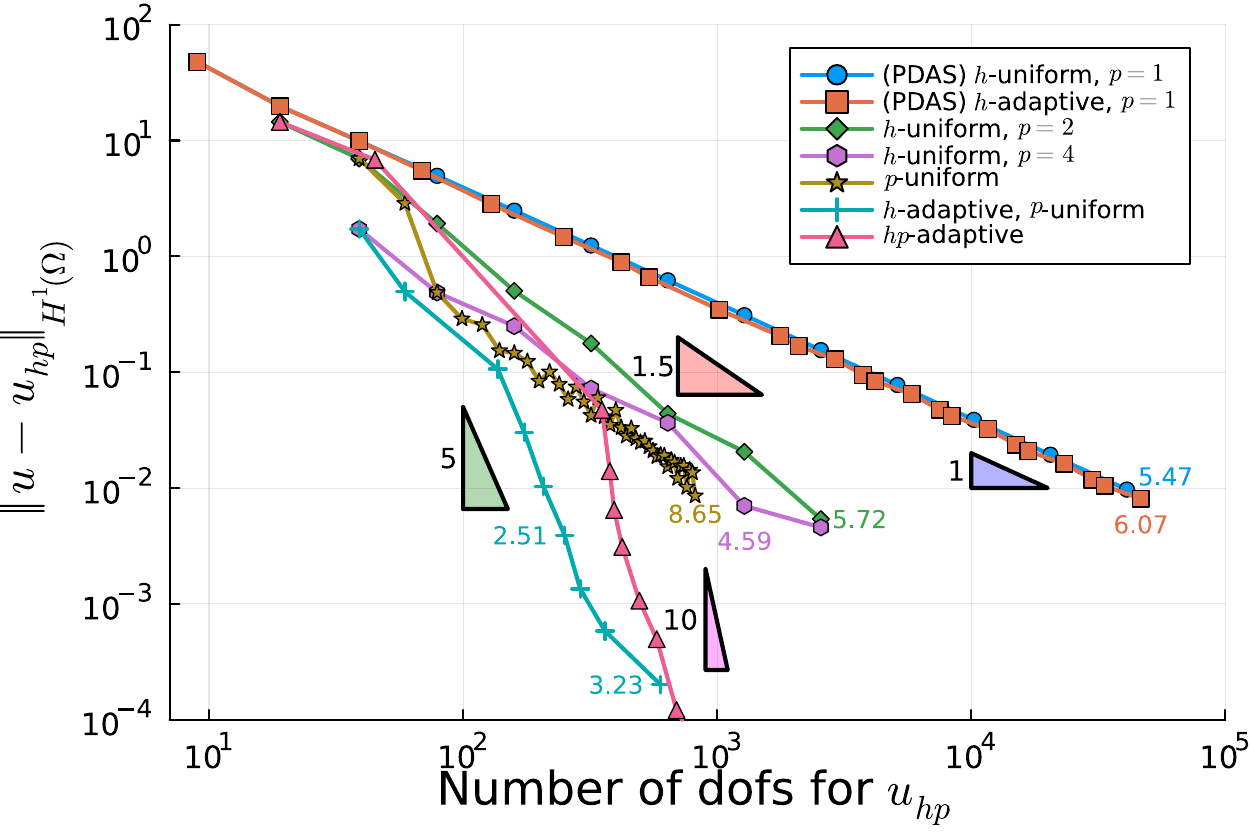}
\caption{(Left) The solution to the 1D obstacle problem with the setup \cref{eq:osc-data-setup}. (Right) Convergence of 7 refinement strategies. $hp$-adaptive and $h$-adaptive follow the strategy in \cite[Sec.~4]{banz2015}, $h$-uniform implies the mesh size is uniformly halved with each refinement and $p$-uniform is where the polynomial degree is incremented by one with each refinement. Any number associated with a data point is the average time taken per Newton iteration, measured in milliseconds, via a direct sparse factorization. In the order listed in the legend, the first two strategies converge at a rate of 1, the next 3 at a rate slightly faster than 3/2 and the final two at a rate of 5 and 10, respectively (in the regime of $h$ and $p$ considered). }
\label{fig:oscillatory-data}
\end{figure}

The first example is a one-dimensional obstacle problem with an oscillatory right-hand side. Let $\omega = 10\pi$ and $c = 2\omega^2$ and consider the parameters:
\begin{align}
\Omega = (0,1), \quad  f(x) = c\sin(\omega x),  \quad \varphi \equiv 1.
\label{eq:osc-data-setup}
\end{align}
The goal is to compare the hpG solver with the primal-dual active set strategy (PDAS) \cite{hintermuller2002} and investigate the effectiveness of an adaptive $hp$-refinement. When using PDAS, $u$ is discretized with continuous piecewise bilinear finite elements on quadrilateral cells. The PDAS solver has no assembly cost at each nonlinear iteration and is effectively solving a Poisson equation over the inactive set of the domain at each iteration. However, each time the active set changes, the stiffness matrix must modified and re-factorized via a sparse Cholesky factorization. For the hpG solver, we use the $\alpha$-update rule $\alpha_1 = 2^{-7}$, $\alpha_{k+1} = \min(\sqrt{2} \alpha_k, 2^{-3})$ and terminate once $\alpha_k = \alpha_{k-1} = 2^{-3}$. We fix the stabilization parameter $\beta = 10^{-8}$ in \cref{eq:pG-matrix}. The results are visualized in \cref{fig:oscillatory-data}.

The strategy labelled ``$hp$-adaptive" in \cref{fig:oscillatory-data} leverages an $hp$-refinement strategy following the  $hp$-a posteriori error estimators developed by Banz and Schr\"oder in \cite[Sec.~4]{banz2015}. The decision on whether to refine a cell or increase the approximation degree is based on the work of Houston and S\"uli \cite[Alg.~1]{houston2005}. As shown, this strategy is very effective at error reduction. In fact, we discover a convergence rate of $O((h/p)^{10})$ for the regimes of $h$ and $p$ that we consider. At the final data point, we have an $H^1$-norm error of $5.61 \times 10^{-5}$, a minimum and maximum discretization degree of 13 and 18 and mesh sizes $h_{\text{min}} = 3.9\times 10^{-4}$ and $h_{\text{max}} = 5 \times 10^{-2}$ across all the cells of the mesh, respectively. Moreover, the hpG iterations are bounded for all ten data points where the solves required 32, 33, 28, 30, 30, 32, 32, 32, 33, and 33 cumulative Newton iterations for each successive refinement, respectively. Recall that on each mesh, the initial guess is initialized as the zero function. Hence, these Newton iterations are \emph{not} dependent on grid-sequencing from the discretized solution of the parent coarser mesh.

The strategy labelled ``$h$-adaptive, $p$-uniform" in \cref{fig:oscillatory-data} also provided a fast convergence rate of $O((h/p)^{5})$. Here we use an adaptive $h$-refinement \cite[Sec.~4]{banz2015}. but coupled with an aggressive $p$-refinement where on the first mesh $p=4$ and then on each successive mesh we increase the polynomial degree by one across all the cells on the mesh. The main advantage of such a strategy is that the transforms required for the assembly and action of nonlinear block $D_\psi$ are parallelized more easily and lead to fast wall-clock solve times. In fact this strategy provided both the smallest $H^1$-norm error as measured by both the number of dofs and the time cost per linear solve, at least in the regime of $h$ and $p$ that we consider. The hpG solves were mesh independent requiring 22 Newton iterations on each refinement. The solve for the finest discretization took 0.087 seconds to achieve an $H^1$-norm error of $2.0 \times 10^{-4}$. Most notably, this strategy required half the average wall-clock time per Newton iteration, when compared to the PDAS strategy, to achieve error smaller than $10^{-2}$. 

The three strategies of uniform $p$-refinement (with a fixed $h$) and uniform $h$-refinement with either $p=2$ or $p=4$ all delivered a convergence rate of slightly faster than $O((h/p)^{3/2})$. This is expected behaviour due to the $H^s(\Omega)$, $s<5/2$, regularity of the solution. In the $p$-refinement, the largest $p$ considered is $p=41$. As expected, both a uniform or adaptive $h$-refinement with $p=1$ via the PDAS solver was capped at a convergence of $O(h)$.

\begin{remark}[Comparing solver times]
We solve the linear systems in the PDAS and hpG strategies via sparse direct solvers. The linear solves in the PDAS could be greatly accelerated, for example, by using algebraic multigrid techniques. However, we contend that a similar improvement is possible for the linear solves in the hpG solver via the preconditioning strategy of \Cref{sec:preconditioning} and other similar strategies \cite{Babuska1991,beuchler2006b,korneev1999,Schwab1998}. To conclude, we are not claiming that a high-order discretization is always the most computationally efficient strategy in terms of error reduction. However,  we have clear evidence to assert that they are certainly competitive.
\end{remark}

\subsection{Oscillatory obstacle}
\label{sec:examples:2d:1}
\begin{figure}[h!]
\centering
\includegraphics[width =0.49 \textwidth]{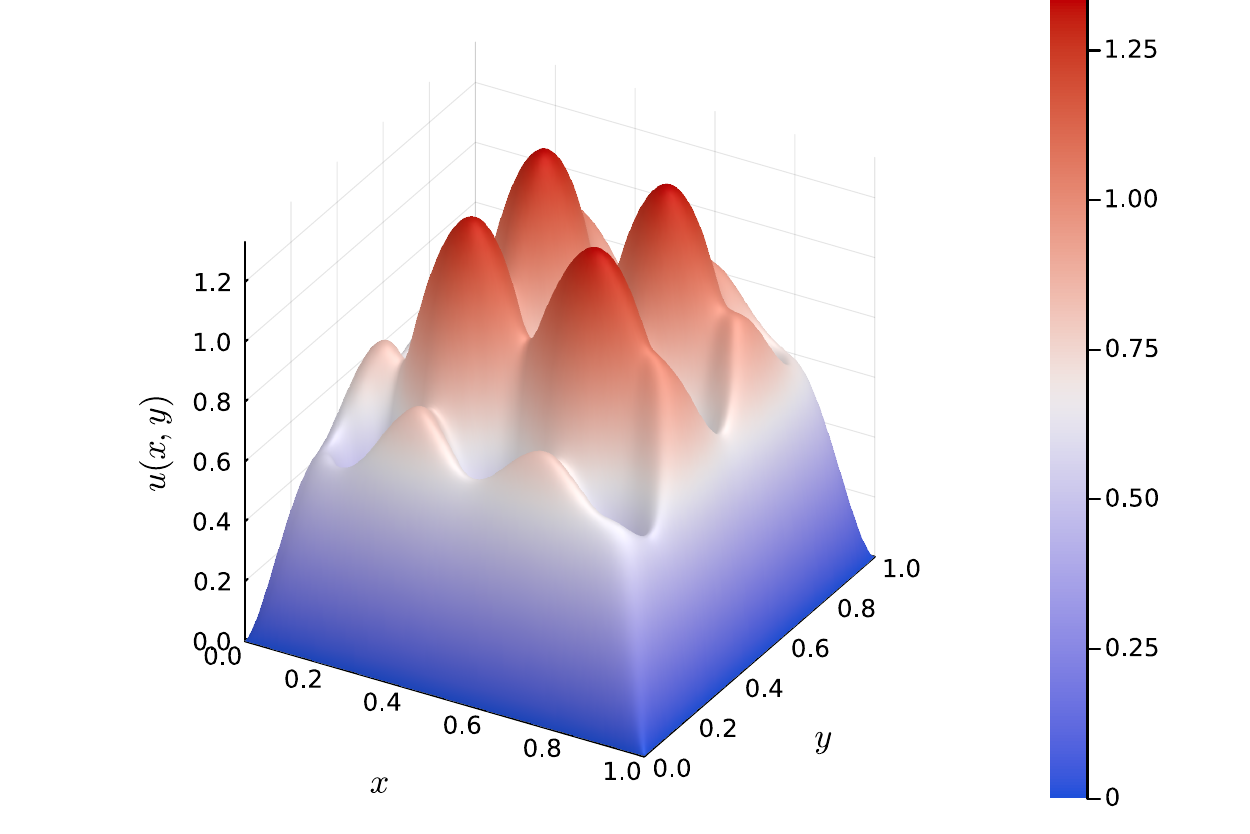}
\includegraphics[width =0.49 \textwidth]{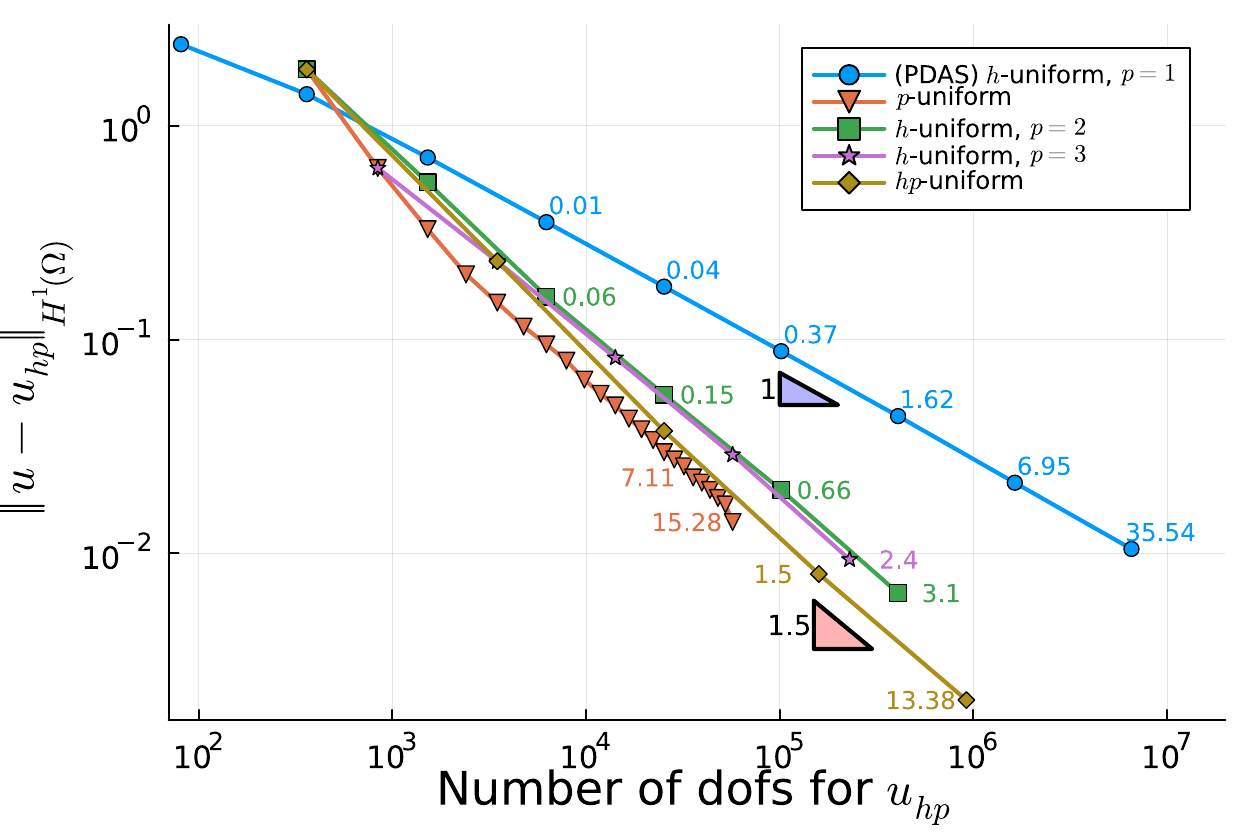}
\caption{(Left) The solution to the obstacle problem with the setup \cref{eq:oscillatory-obstacle-setup}. (Right) Approximate $H^1$-norm error of five strategies for the obstacle problem with the setup \cref{eq:oscillatory-obstacle-setup}. All hpG solves required a total of 24 Newton iterations independently of $p$ and $h$. The numbers attached to a data point are the average time taken per linear solve, as measured in seconds, via a sparse factorization. The higher order discretizations achieve a smaller error with fewer dofs and a faster solve time. The highest partial degree plotted is $p=24$. The triangles indicate rates of convergence. When $p=1$ we observe an $O(h)$ rate of convergence whereas for the other strategies the rate is roughly $O((h/p)^{3/2})$.}
\label{fig:oscillatory-obstacle}
\end{figure}

Consider the following setup of a two-dimensional obstacle problem:
\begin{align}
\Omega = (0,1)^2, \;\;  f(x,y) = 100, \;\; \text{and} \;\;
\varphi(x,y) =(1+J_0(20x))(1+J_0(20y)),
\label{eq:oscillatory-obstacle-setup}
\end{align}
where $J_0$ denotes the zeroth order Bessel function of the first kind \cite[Sec.~10.2(ii)]{NIST:DLMF}. We are not aware of the exact solution of this problem and, therefore, estimate the error against a heavily-refined discretization that is plotted in \cref{fig:oscillatory-obstacle}. In this example we focus primarily on uniform $h$ and $p$-refinements and discover that the best ratio of (approximate) error to computational expense is achieved by a $hp$-uniform refinement, i.e.~with each refinement the mesh size is halved and the partial discretization degree is incremented by one. We believe that this subsection provides a clear counterexample to the notion that low-order discretizations should always be preferred.

We initialize the discretization by meshing the domain into $10\times 10$ uniform quadrilateral cells and use the same $\alpha$-update rule as in \Cref{sec:examples:1d}. We do not use any stabilization, i.e.~$\beta = 0$ in \cref{eq:pG-matrix}. We find that all hpG runs required 24 Newton iterations independently of $h$ and $p$. We plot the (approximate) $H^1$-norm error of the three strategies in \cref{fig:oscillatory-obstacle} and include the average time taken per linear solve for each solver (including any matrix assembly costs).

The convergence of the PDAS solver is capped at $O(h)$ since the partial degree is fixed at $p=1$ whereas the other strategies observe an expected convergence rate of $O((h/p)^{3/2})$. Per dof, the high-order discretizations achieve the smallest error. Moreover, the average time taken per linear solve is clearly competitive with the low-order discretizations coupled with a finer mesh. In fact the best strategy, $hp$-uniform refinement, featured a discretization that resulted in wall-clock times that were roughly 24 times faster, per linear solve, than the PDAS strategy for an error of $10^{-2}$.

We now test our linear system solver strategy via the Schur decomposition \cref{fig:solver2}. In \cref{tab:oscillatory-obstacle} we tabulate the average right-preconditioned GMRES iterations and wall-clock time per Newton step of the hpG solver for various choices of $h$ and $p$. We pick $\beta = 10^{-4}$ in $E_\beta$ in \cref{eq:pG-matrix} and choose a relative stopping tolerance of $10^{-5}$ for the GMRES solver. For comparison, we include a row with the average wall-clock time per Newton step via a sparse LU factorization.  The GMRES iterations are bounded with respect to $p$ and only grow at a mild polylogarithmic rate as $h \to 0$. Moreover, the GMRES solver always resulted in faster average solve times per linear iteration than the LU factorization.

\begin{table}[ht]
\small
\centering
\scalebox{0.7}{
\begin{tabular}{|c|cccc|ccc|ccc|ccc|}
\hhline{~|----|---|---|---|}
\rowcolor{lightgray!10} \multicolumn{1}{c|}{\cellcolor{lightgray!00}} &\multicolumn{4}{c|}{$p=n$, $h=1/10$} & \multicolumn{3}{c|}{$p=2$, $h=2^{-n}/10$}  & \multicolumn{3}{c|}{$p=3$, $h=2^{-n}/10$}  & \multicolumn{3}{c|}{$p=n+2$,  $h = 2^{-n}/10$}\\
\hline
\cellcolor{lightgray!10} $n$ & 3 &10  & 17 & 24   & 1& 3& 5& 0& 2  & 4 & 0& 2 & 4  \\ 
\hline
\cellcolor{lightgray!10} GMRES Its.&  16.38&23.63 &23.54&23.38&  20.33 & 32.08 & 36.54 & 16.38& 26.79 &35.13&14.21&29.33&35.29\\
\hline
\cellcolor{lightgray!10} GMRES Time &  0.00 &0.19 &1.67&7.99&  0.00 & 0.12 & 2.64 & 0.00& 0.06 &1.58&0.00&0.23&12.07\\
\hline
\cellcolor{lightgray!10} LU Time &  0.00 & 0.22 &2.75 &15.28&  0.01 & 0.15 & 3.10 &  0.01& 0.10 &2.40&0.00&0.20&13.38\\
\hline
\end{tabular}}
\caption{A comparison of the average right-preconditioned GMRES iterations and wall-clock timings (in seconds) per Newton iteration for various choices of $h$ and $p$ over a run of the hpG algorithm to solve the obstacle problem with the setup \cref{eq:oscillatory-obstacle-setup}. We utilize the preconditioner  outlined in \cref{fig:solver2} with a GMRES relative tolerance of $10^{-5}$. We use the stabilization $E_\beta$ with $\beta = 10^{-4}$ and include the average wall-clock timings of a sparse LU factorization for comparison. We observe a bounded iteration count with respect to $p$ and a mild polylogarithmic growth with respect to $h$. The GMRES solver is almost always faster than the LU factorization.} 
\label{tab:oscillatory-obstacle}
\end{table}

\subsection{A gradient-type constrained problem}
\label{sec:examples:gradient}
Consider the gradient-type constrained Dirichlet minimization problem in \cref{eq:ad:gradient} and fix the parameters as
\begin{align}
\Omega = (0,1)^2, \;\;  f(x,y) = 20, \;\; \text{and} \;\;
\varphi(x,y) =
\begin{cases}
1/2 & \text{if} \; x \; \text{or} \; y \in [0,1/4] \cup [3/4,1], \\
\infty  &\text{otherwise}.
\end{cases}
\label{eq:gradient-setup}
\end{align}

\begin{figure}[h!]
\centering
\includegraphics[width =0.45 \textwidth]{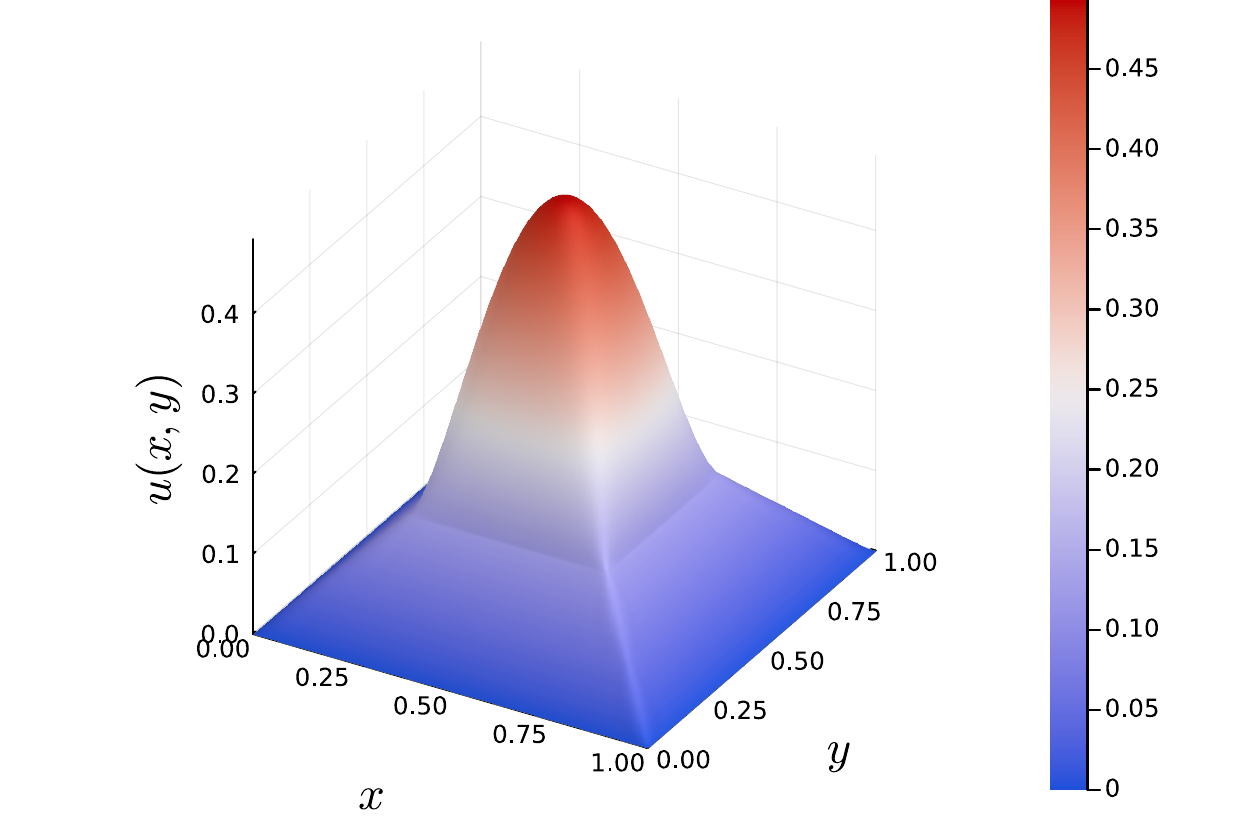}
\includegraphics[width =0.45 \textwidth]{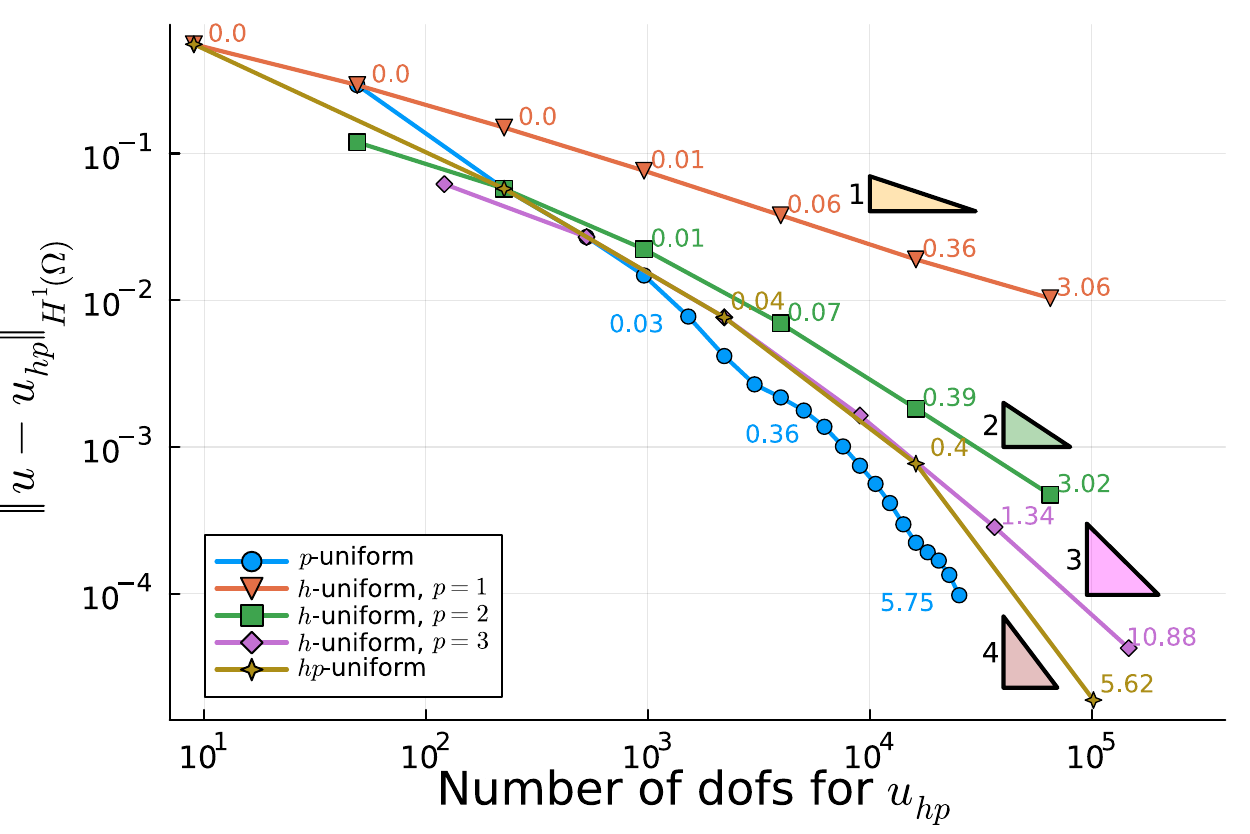}
\caption{(Left) The solution to the gradient-constraint problem \cref{eq:ad:gradient} with the setup \cref{eq:gradient-setup}. (Right) Approximate $H^1$-norm error of five strategies. All hpG solves required between 74 and 78 Newton iterations. The numbers attached to data points are the average time taken per linear solve, as measured in seconds, with a sparse LU factorization. The highest partial degree plotted is $p=20$. The triangles indicate rates of convergence. We observe $O(h^p)$ convergence in the range of $h$ and $p$ considered.}
\label{fig:gradient-constraint}
\end{figure}

In \cref{fig:gradient-constraint} we plot the solution of the gradient-constrained problem with the setup \cref{eq:gradient-setup}. There is a kink in the solution at $x , y \in \{1/4, 3/4\}$ which coincides directly with the definition of $\varphi$. Once again, we are not aware of the exact solution of this problem and, therefore, we measure the convergence against a heavily-refined discretization. We initialize the discretization by meshing the domain into $4\times 4$ uniform quadrilateral cells (or $8 \times 8$ in the $p$-uniform strategy) which aligns with the discontinuity of $\varphi$. We use the $\alpha$-update rule $\alpha_1 = 2^{-7}$, $\alpha_{k+1} = \min(\sqrt{2} \alpha_k, 2^{2})$ and terminate once $\alpha_k = \alpha_{k-1} = 2^{2}$. The hpG solver features $hp$-robustness with all hpG solves requiring between 74 and 78 Newton iterations irrespective of $p$ and $h$. The convergence plot is also provided in \cref{fig:gradient-constraint}. Surprisingly, in the strategies that also refine the mesh, we observe regions with $O(h^p)$ convergence rates in contrast to the obstacle problems considered previously. We hypothesize that this superconvergence arises from the alignment of the mesh with the discontinuity in $\varphi$. In this example, this alignment coincides with the transition between the active and inactive regions of the solution which is where the loss of regularity occurs. Hence, the solution restricted to each cell $K$ is more regular than $H^{5/2}(K)$ and, therefore, benefits from the high-order discretization. Once again, we observe significant speedups when utilizing discretizations $p \geq 2$. We observe that the induced linear systems are around 100 times faster to solve for the strategies where $p > 1$ in order to reach an error of $10^{-2}$.

In \cref{tab:gradient-constraint} we report the effectiveness of the Schur decomposition \cref{fig:solver2} for solving the Newton systems. The average GMRES iterations per linear solve are $p$-robust and grow with up to a linear rate with $1/h$. The growth is exacerbated by the choice of an adaptive stabilization parameter where $\beta = 10^{-p+\log_2 h} \to 0$  as $p \to \infty$ and $h \to 0$. If $\beta$ is fixed to be constant, then the growth is at most polylogarithmic but causes the error to plateau. For all reported choices of $h$ and $p$, we observe that the average wall-clock timing per linear solve of the preconditioned GMRES strategy is faster than a direct sparse LU factorization.

\begin{table}[h!]
\small
\centering
\scalebox{0.7}{
\begin{tabular}{|c|cccc|ccc|ccc|cc|}
\hhline{~|----|---|---|--|}
\rowcolor{lightgray!10} \multicolumn{1}{c|}{\cellcolor{lightgray!00}} &\multicolumn{4}{c|}{$p=n$, $h=2^{-3}$} & \multicolumn{3}{c|}{$p=1$, $h=2^{-n}$}  & \multicolumn{3}{c|}{$p=3$, $h=2^{-n}$}  & \multicolumn{2}{|c|}{$p=n-1$,  $h = 2^{-n}$}\\
\hline
\cellcolor{lightgray!10} $n$ & 5 &10  & 15 & 20   & 4& 6& 8& 3& 5  & 7 &  4 & 6  \\ 
\hline
\cellcolor{lightgray!10} GMRES Its.&  61.76 &71.76 &72.03&72.00& 10.31 & 20.73 & 33.73 & 21.26 & 50.71  &72.76&34.73&77.48\\
\hline
\cellcolor{lightgray!10} GMRES Time & 0.03 &  0.33 & 1.70 & 5.45 &  0.00 & 0.02 & 0.83 & 0.00& 0.16  & 4.49&0.03&4.66\\
\hline
\cellcolor{lightgray!10} LU Time &  0.03 & 0.36 &2.60 &5.75& 0.00  & 0.06  & 3.06  & 0.01 & 0.19  & 10.88 &0.04& 5.62\\
\hline
\end{tabular}}
\caption{A comparison of the average right-preconditioned GMRES iterations and wall-clock timings (in seconds) per Newton iteration for various choices of $h$ and $p$ over a run of the hpG algorithm to solve \cref{eq:ad:gradient} with the setup \cref{eq:gradient-setup}. We utilize the preconditioning strategy outlined in \cref{fig:solver2} with a GMRES relative tolerance of $10^{-3}$ and terminate GMRES if it reaches 150 iterations. We use the stabilization $E_\beta$ with $\beta = 10^{-p + \log_2 h }$ and include the average wall-clock timings obtained with a sparse LU factorization for comparison. We observe that the iterations are $p$-robust and grow with up to a mild linear rate with $1/h$. The preconditioning strategy always offers a speedup when compared to a direct sparse LU factorization.} 
\label{tab:gradient-constraint}
\end{table}

\subsection{An obstacle-type quasi-variational inequality}
\label{sec:examples:qvi}

In this example we consider the thermoforming \emph{quasi-variational} inequality (QVI). In a QVI, the obstacle $\varphi$ is dependent on the solution itself. Solvers for QVIs posed in an infinite-dimensional setting are scarce and most examples in the literature are solved by means of a fixed point iteration, penalty or augmented Lagrangian technique \cite{kanzow2019,AHR}. Very recently a semismooth Newton method was introduced for a class of obstacle-type QVIs in \cite{alphonse2024}. The solver in \cite{alphonse2024} requires the realization of an active set and, therefore, is restricted to low-order FEM discretizations. The QVI can also be directly tackled by an extension of the pG algorithm \cite[Sec.~3.5]{dokken2024} but handling the nonlinear terms such that the discretizations remain sparse as $p \to \infty$ is nontrivial. Hence in this example we opt for a fixed point approach where the obstacle subproblems are solved via the hpG solver. We believe this is the highest order discretization of an elliptic obstacle-type QVI that is reported in the literature.

Given a $\Phi_0 \in H^1(\Omega)$, $\xi \in C^2(\bar \Omega) \cap H^1_0(\Omega)$, $f \in L^2(\Omega)$, $\gamma > 0$, and a globally Lipschitz and nonincreasing function $g : H^1(\Omega) \to L^2(\Omega)$, the thermoforming problem is given by
\begin{equation}
\label{eq:modelqvi}
\begin{aligned}
	&\text{Find $u \in H_0^1(\Omega)$  satisfying, for all $v \in \{ w \in H_0^1(\Omega) : w \leq \Phi(u) \coloneqq \Phi_0 + \xi T\}$,} 
\\
&\qquad\qquad\quad u \leq  \Phi(u),
\quad(\nabla u, \nabla (v-u))_{L^2(\Omega)} \geq ( f, v-u)_{L^2(\Omega)},
\\
&\text{with $T \in H^1(\Omega)$ satisfying, for all $q \in H^1(\Omega)$,}
\\
&\qquad\qquad\quad  (\nabla T, \nabla q)_{L^2(\Omega)} + \gamma ( T, q)_{L^2(\Omega)} = (g(\Phi_0 + \xi T -u), q)_{L^2(\Omega)}.
\end{aligned}
\end{equation}

\begin{remark}[Thermoforming QVI]
In two dimensions, $\Omega \subset \mathbb{R}^2$, \eqref{eq:modelqvi}
provides a simple model for the problem 
of determining the displacement $u \in H^1_0(\Omega)$ of 
an elastic membrane,
clamped at the boundary $\partial \Omega$,
that has been heated, and is 
pushed by means of an external force $f \in L^2(\Omega)$
into a metallic mould with original shape $\Phi_0 \in H^1(\Omega)$ and final shape $\Phi(u) \in H^1(\Omega)$. The deformation is due to the mould's temperature field $T \in H^1(\Omega)$ which varies according to the
membrane's temperature. The heat transfer is modelled by a conduction coefficient $\gamma > 0$, a given globally Lipschitz and nonincreasing function $g :H^1(\Omega) \to L^2(\Omega)$, and a smoothing function $\xi \in H^1_0(\Omega) \cap C^2(\bar \Omega)$ that incorporates the distance
between the membrane and the mould. For more details on the thermoforming problem, including 
its derivation and its background,  
we refer to \cite[Section 6]{AHR} and the references therein.
\end{remark}

The fixed point approach proceeds as follows, for $T_0 \equiv 0$ and $i\in \mathbb{N}_0$, repeat the following two steps until convergence:
\begin{enumerate}[label=(\Roman*)]
\item \label{qvi:step1} Given a $T_i \in H^1(\Omega)$, solve the obstacle problem, find $u_{i+1} \in H_0^1(\Omega)$  satisfying  $\forall v \in \{w \in H_0^1(\Omega) : w \leq \Phi_0 + \xi T_i \; \text{a.e.}\}$,
\begin{align}
u_{i+1} \leq \Phi_0 + \xi T_i,
\; (\nabla u_{i+1}, \nabla (v-u_{i+1}))_{L^2(\Omega)} \geq ( f, v-u_{i+1})_{L^2(\Omega)}.
\end{align}
\item \label{qvi:step2} For all $q \in H^1(\Omega)$, solve the (nonlinear) PDE for $T_{i+1} \in H^1(\Omega)$:
\begin{align}
(\nabla T_{i+1}, \nabla q)_{L^2(\Omega)} + \gamma ( T_{i+1}, q)_{L^2(\Omega)} = (g(\Phi_0 + \xi T_{i+1} -u_{i+1}), q)_{L^2(\Omega)}.
\end{align}
\end{enumerate}

\begin{figure}[h!]
\centering
\includegraphics[width =0.34 \textwidth]{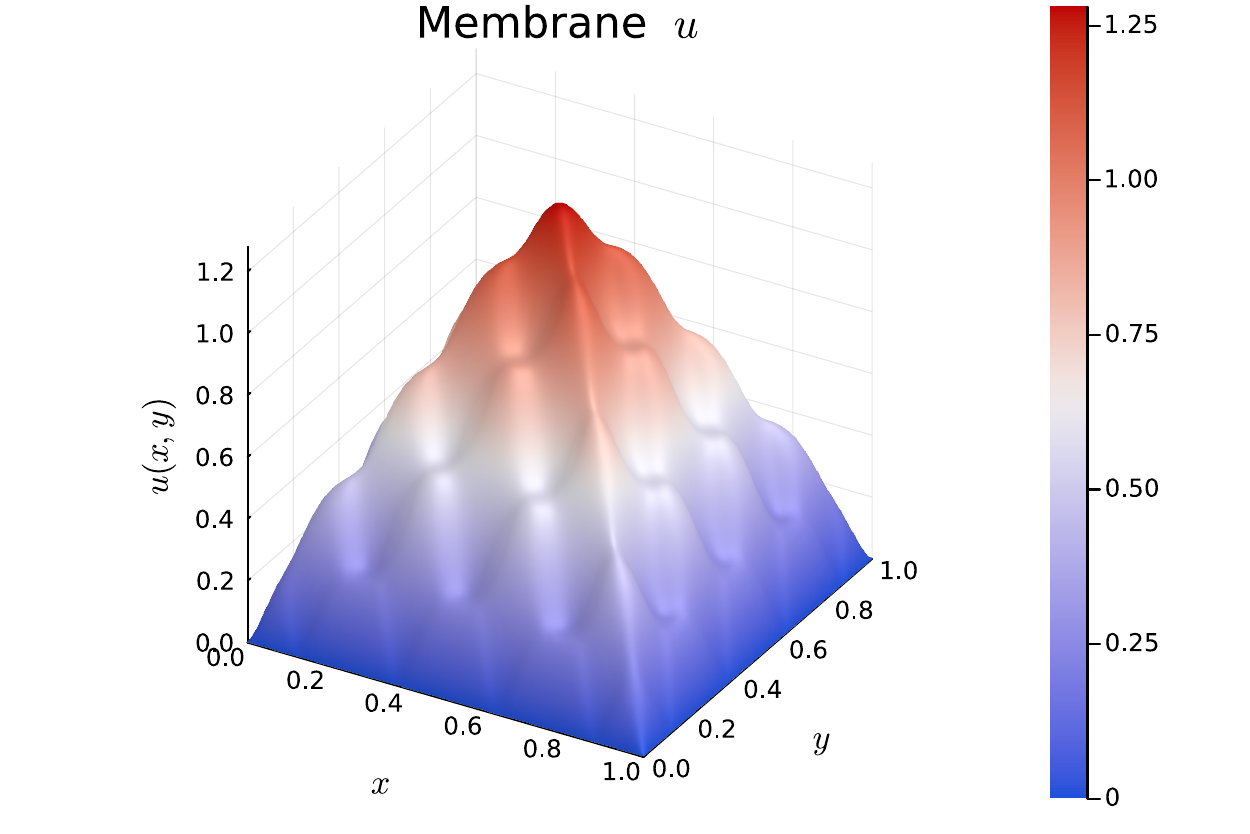}
\includegraphics[width =0.34 \textwidth]{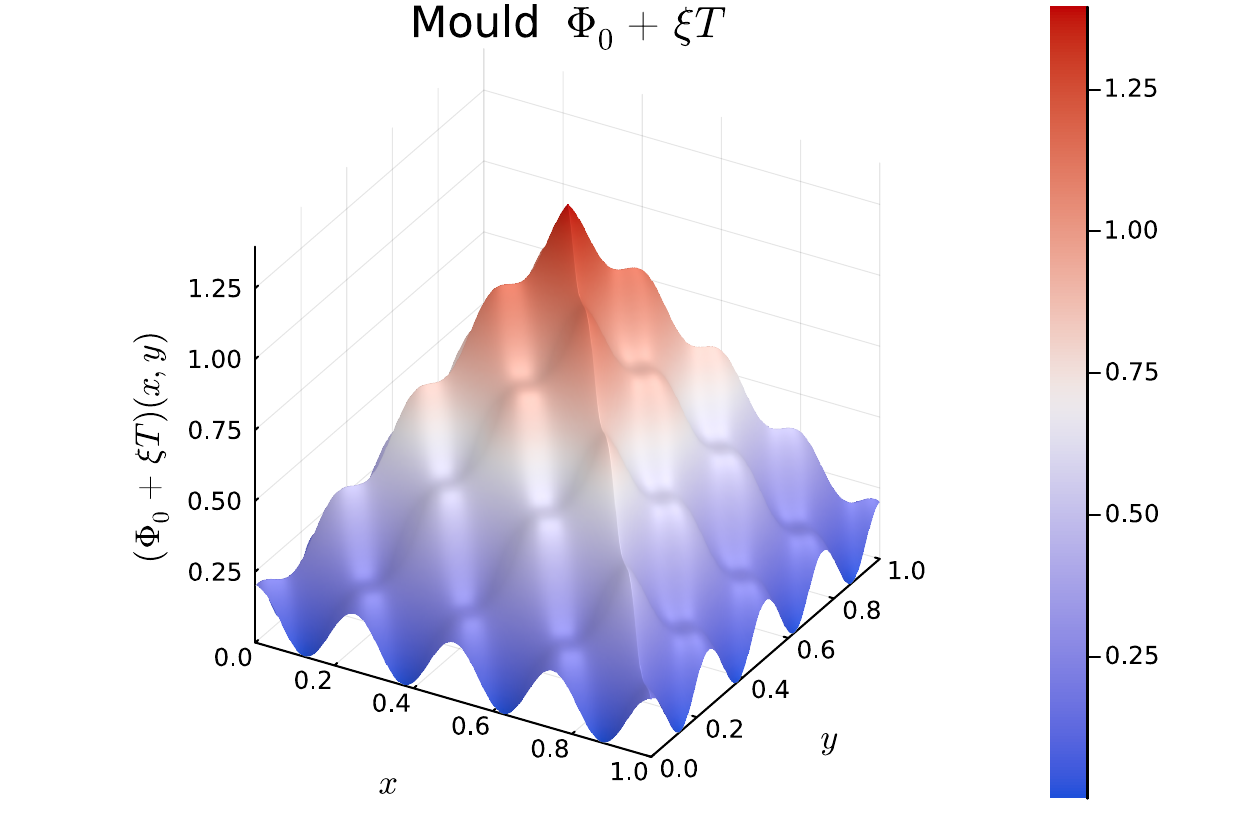}
\includegraphics[width =0.3 \textwidth]{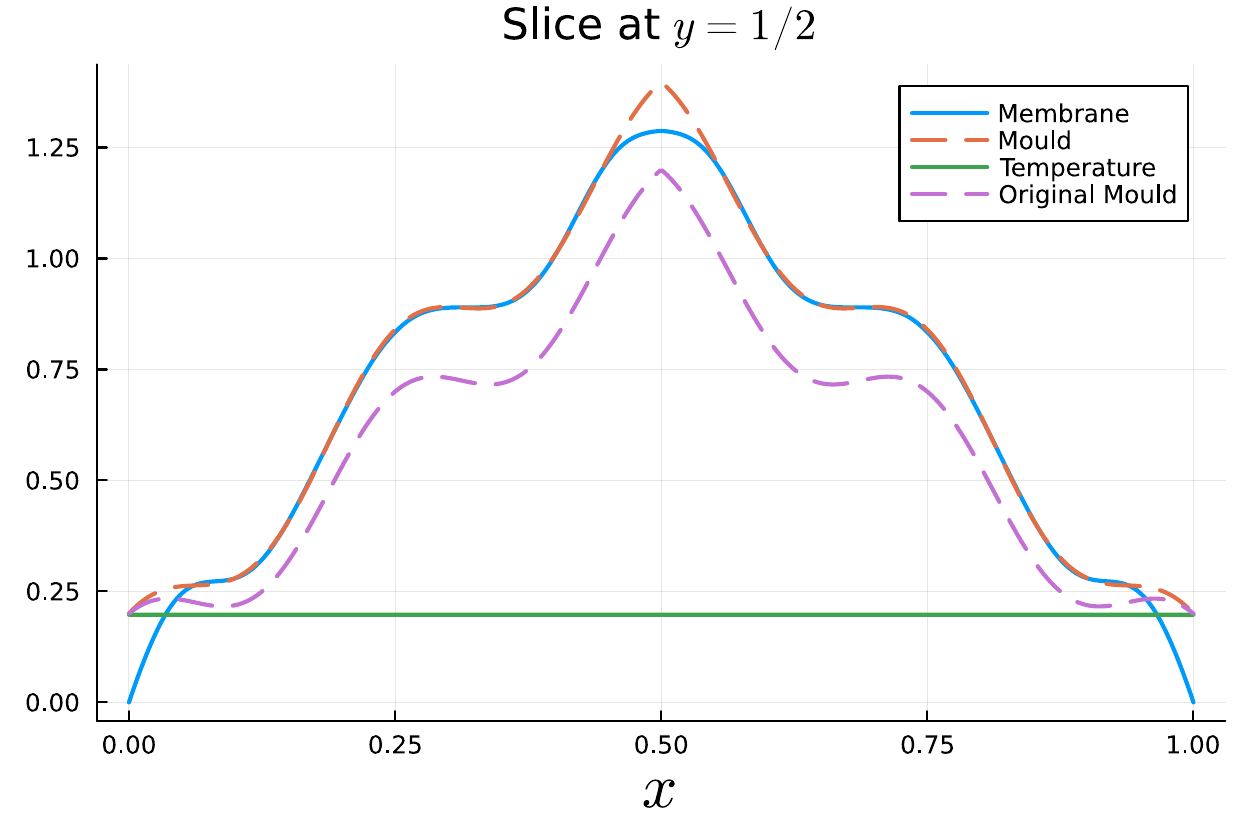}
\caption{Surface plots of the membrane $u$ (left) and the final mould (middle) of the thermoforming problem \cref{eq:modelqvi} with setup \cref{eq:qvi-setup}. On the right is a slice plot at $y=1/2$ of the membrane and final mould as well as the original mould $\Phi_0$ and final temperature $T$.}
\label{fig:thermoforming}
\end{figure}

We now fix the example parameters as:
\begin{equation}
\label{eq:qvi-setup}
\begin{gathered}
\Omega  = (0,1)^2,
\qquad
f(x,y) = 100,  \qquad \xi(x,y)
=
\sin(\pi x)
\sin(\pi y)
,\qquad
\gamma =1,\\
\Phi_0(x, y)
=
11/10 - 2 \max(|x - 1/2|,|y - 1/2|) + \cos(8\pi x)\cos(8 \pi y)/10,
\\
g(s) = \begin{cases}
1/5 & \text{if} \; s\leq 0,\\
(1-s)/5& \text{if} \; 0< s < 1,\\
0 & \text{otherwise}.
\end{cases}
\end{gathered}
\end{equation}
The data \cref{eq:qvi-setup} is such that we are guaranteed the existence and uniqueness of a solution to \cref{eq:modelqvi} \cite{AHR}. Moreover, one can show there is a global contraction in the fixed point iteration induced by steps \labelcref{qvi:step1,qvi:step2}, cf.~\cite[Sec.~4.4]{alphonse2024} and, therefore, the scheme is guaranteed to converge to the (unique) solution. We plot the resulting membrane $u$ and mould $\Phi_0 + \xi T$ in \cref{fig:thermoforming} as well as a slice at $y=1/2$.

We mesh the domain into a $4 \times 4$ uniform mesh (i.e.~16 cells total) and discretize $(u, \psi, T)$ with $(u_\hp, \psi_\hp, T_\hp) \in U_{h, \vect{p}_x,\vect{p}_y}  \times \Psi_{h, \vect{p}_x-2,\vect{p}_y-2} \times U_{h, \vect{p}_x,\vect{p}_y}$. We terminate the fixed point algorithm once $\| u_{i} - u_{i-1} \|_{H^1(\Omega)} \leq 3 \times 10^{-3}$.

To solve the obstacle problem in step \labelcref{qvi:step1}, we use the hpG solver with the $\alpha$-update rule $\alpha_1 = 2^{-6}$, $\alpha_{k+1} = 4 \alpha_k$, and terminate once $\alpha_k =1$. We use the stabilization $E_\beta$ in \cref{eq:pG-matrix} with $\beta = 10^{-6}$. To solve the hpG Newton linear systems we  leverage the Schur decomposition in \cref{fig:solver2}. We choose a GMRES relative and absolute  stopping tolerance of $10^{-7}$ for approximating the inverse of $S$.

We found that an efficient nonlinear solver for step \labelcref{qvi:step2} was  Newton's method where the linear systems are solved by GMRES left-preconditioned with $A + \gamma M_u$ where $A$ and $M_u$ are the stiffness and mass matrices for $U_\hp$. Hence, the assembly of the (increasingly dense as $p \to \infty$) Jacobian is never required and the Jacobian-vector product can be computed in quasi-optimal complexity (similar to the discussion in \Cref{sec:implementation}). Since $A + \gamma M_u$ is very sparse, it can be factorized efficiently with a sparse Cholesky factorization (the factorization took 0.2 seconds with $p=82$). We chose a relative stopping tolerance of $1.5 \times 10^{-8}$ for the GMRES solver.%

In \cref{tab:thermoforming} we provide the iteration counts of the outer fixed point method as well as the average number of Newton iterations per fixed point iteration and the average number of preconditioned GMRES iterations to invert $S$ per Newton iteration for  steps \labelcref{qvi:step1,qvi:step2} with increasing partial degree $p$. We observe $p$-independent Newton iteration counts for both steps \labelcref{qvi:step1,qvi:step2}. The average number of GMRES iterations per Newton iteration appears to be bounded above by 22.4 and 3.11 for  steps \labelcref{qvi:step1,qvi:step2}, respectively.

\begin{table}[h!]
\small
\centering
\scalebox{0.9}{
\begin{tabular}{|c|c|c|c|c|c|}
\hhline{~~|--|--|}
\rowcolor{lightgray!10} \multicolumn{1}{c}{\cellcolor{lightgray!00}}     &   \multicolumn{1}{c|}{\cellcolor{lightgray!00}}  &\multicolumn{2}{c|}{Step \labelcref{qvi:step1}}  & \multicolumn{2}{c|}{Step  \labelcref{qvi:step2}} \\ 
\hline
\rowcolor{lightgray!10} $p$  & Fixed point & Avg.~Newton & Avg.~GMRES & Avg.~Newton & Avg.~GMRES\\
\hline 
6 &   4&   15.00 & 11.00 & 1.50 & 2.83 \\
12 & 4 & 15.25 &15.85 & 2.00& 3.13 \\
22 & 4 & 16.00 & 19.36 & 2.00 &3.00 \\
32 & 4 & 16.00& 21.09 & 2.00 &3.00 \\
42 & 4 & 15.75 &21.75 & 2.25 &3.11 \\
52 & 4 & 15.00 &22.40 & 2.00 &3.00 \\
62 & 4 & 15.00 & 21.90 & 2.00 & 3.00\\
72 & 4 & 15.00 & 21.90 & 2.00 & 3.00\\
82 & 4 & 15.25 & 21.61 & 2.00 & 3.00\\
\hline
\end{tabular}}
\caption{The partial degree $p$, the number of outer iterations of the fixed point scheme as well as the average number of Newton iterations per fixed point iteration and average number of preconditioned GMRES iterations to invert $S$ per Newton iteration in steps \labelcref{qvi:step1,qvi:step2} to approximate the solution of the QVI \cref{eq:modelqvi} with setup \cref{eq:qvi-setup}. The algorithm terminates once $\| u_{i} - u_{i-1} \|_{H^1(\Omega)} \leq 3 \times 10^{-3}$. The outer fixed point loop is $p$-independent and we observe bounded iteration counts in the Newton and GMRES solvers for both steps \labelcref{qvi:step1,qvi:step2}.} 
\label{tab:thermoforming}
\end{table}

\subsection{A three-dimensional obstacle problem}
\label{sec:examples:3d}

As evidenced in the previous examples, smaller errors are achieved with significantly fewer dofs with the hpG solver than with a $p=1$ discretization. Moreover, the solve times were also faster, including any assembly costs. In this example, we examine the costs of the hpG solver applied to a three-dimensional obstacle problem. 

As an LU factorization of the Newton linear systems may quickly become infeasible, we immediately focus on the iterative solvers developed in \Cref{sec:preconditioning}. The first concern is that the one-time Cholesky factorization of the stiffness matrix may become infeasible on fine meshes and higher degrees. Thus, in \cref{tab:cg}, we investigate the effectiveness of inverting $A_\alpha$ via CG preconditioned with a Ruge--St\"uben algebraic multigrid cycle \cite[Ch.~4]{ruge1987}. We observe that the Runge--St\"uben AMG setup is magnitudes faster than a Cholesky factorization but, as expected, applying preconditioned CG is slower than leveraging the Cholesky factorization to apply the inverse.  The CG iteration counts undergo an initial growth with both $h$ and $p$-refinement, although the growth appears to stagnate. For scenarios where the initial growth is problematic, we recommend considering $p$-multigrid which will likely produce more robust iteration counts \cite{Brubeck2022}.
\begin{table}[h!]
\small
\centering
\scalebox{0.9}{
\begin{tabular}{|c|ccc|ccc|}
\hhline{~|---|---|}
\rowcolor{lightgray!10} \multicolumn{1}{c|}{\cellcolor{lightgray!00}} &\multicolumn{3}{c|}{Cholesky factorization} & \multicolumn{3}{c|}{Apply Cholesky inverse} \\
\hline
\rowcolor{lightgray!10} $p \; \backslash$ Cells & $4^3$ &$8^3$  & $16^3$ & $4^3$ &$8^3$  & $16^3$  \\ 
\hline
\cellcolor{lightgray!10} 2& $1.23 \text{e-3}$ & $1.49 \text{e-2}$ &0.637 & $2 \text{e-5}$&   $3.10 \text{e-4}$ &  $1.62 \text{e-2}$\\
\hline
\cellcolor{lightgray!10} 3 & $4.57 \text{e-3}$  &$0.215$  & 17.6 & $8 \text{e-5}$& $4.99 \text{e-3}$  &0.231\\
\hline
\cellcolor{lightgray!10} 4&  $1.29 \text{e-2}$  &$0.994$ & 126.1& $2.10 \text{e-4}$& $1.63 \text{e-2}$ & 0.718\\
\hline
\cellcolor{lightgray!10} 5&  $3.86 \text{e-2}$  &6.71 & 322.51& $1.95 \text{e-3}$& $4.24 \text{e-2}$ &1.58\\
\hline
\cellcolor{lightgray!10} 6&  $9.02 \text{e-2}$  &37.78& -& $4.18 \text{e-3}$& 0.227 &-\\
\hline
\cellcolor{lightgray!10} 7&  $0.251$  &53.57& -& $5.44 \text{e-3}$& 0.362 &-\\
\hline
\rowcolor{lightgray!10} \multicolumn{1}{c|}{\cellcolor{lightgray!00}} &\multicolumn{3}{c|}{AMG Setup} & \multicolumn{3}{c|}{Preconditioned CG} \\
\hline
\rowcolor{lightgray!10} $p \; \backslash$ Cells & $4^3$ &$8^3$  & $16^3$ & $4^3$ &$8^3$  & $16^3$  \\ 
\hline
\cellcolor{lightgray!10} 2& $7.79 \text{e-4}$ & $1.20\text{e-2}$& $1.09\text{e-1}$  & $4.27 \text{e-3}$ (20)& $2.90 \text{e-2}$ (29)&0.663 (33)\\
\hline
\cellcolor{lightgray!10} 3 &  $5.27 \text{e-3}$  & $6.76 \text{e-2}$ &  $1.030$ &$9.70 \text{e-3}$ (29)&$0.294$ (31)& 2.61 (32)\\
\hline
\cellcolor{lightgray!10} 4&  $8.46 \text{e-3}$& $0.207$ & 2.22  &$5.70 \text{e-2}$ (76)&$1.54$ (96)&$17.89$ (103)\\
\hline
\cellcolor{lightgray!10} 5&  $1.40 \text{e-2}$  &$0.172$ & 3.50& $0.142$ (90)& $3.32$ (98) &$28.92$ (98)\\
\hline
\cellcolor{lightgray!10} 6&  $ 1.77\text{e-2}$  &0.232&3.91 & $0.418$ (165)& $10.23$ (204) &101.46 (216)\\
\hline
\cellcolor{lightgray!10} 7&  $ 2.37\text{e-2}$  &0.375& 6.97& $ 0.870$ (180)& $15.04$ (205)&137.23 (208)\\
\hline
\end{tabular}}
\caption{Wall-clock timings, in seconds, to perform a Cholesky factorization of $A = LL^\top$, apply the inverse to a vector of ones $L^{-\top}L^{-1} \mathbf{1}$, construct the Runge--St\"uben AMG setup on $A$, and compute $A^{-1} \mathbf{1}$ to absolute and relative tolerances of $10^{-6}$ via CG preconditioned with the AMG cycle. The bracketed numbers are the number of preconditioned CG iterations required. A dash indicates that the Cholesky factorization failed due to insufficient memory on a machine with 16GB of RAM.} 
\label{tab:cg}
\end{table}

Since, unlike the two-dimensional case, $A_\alpha$ may present a computational bottleneck, we choose to investigate the cost of alternative iterative solvers from \cref{fig:solver} rather than just \cref{fig:solver2}. The goal is to balance the number of outer FGMRES iterations with the number of $A_\alpha$ inverses dictated by $P = P_F, P_D$ or $P_L P_D$ as well as the number of inner GMRES iterations for $S$, noting that each inner GMRES iteration also requires an $A_\alpha$ inverse.

We choose the obstacle problem setup
\begin{align}
\Omega = (0,1)^3, \;\;  f(x,y,z) = 40, \;\; \text{and} \;\;
\varphi(x,y,z) =1.
\end{align}
We use the $\alpha$-update rule $\alpha_1 = 2^{-2}$, $\alpha_{k+1} = \min(\sqrt{2} \alpha_k, 2)$ and terminate once $\alpha_k = \alpha_{k-1} = 2$ constituting six outer proximal steps. We use the stabilization $E_\beta$ in \cref{eq:pG-matrix} with $\beta = 10^{-5}$ and use an absolute and relative tolerance of $10^{-5}$ for the outermost iterative linear system solver.

For smaller mesh sizes and higher $p$, where a preconditioned CG iterative is required to invert $A_\alpha$, then the fastest strategy is the one that results in the fewest $A_\alpha$-inverses per Newton step. According to \cref{tab:3d}, this appears to be the Schur decomposition \cref{fig:solver2}. It is better to forgo the outer FGMRES loop and prioritise inverting the Schur complement to sufficient accuracy to recover the Newton updates in one step.

\begin{table}[h!]
\small
\centering
\scalebox{0.85}{
\begin{tabular}{|c|c|c|c|c|c|}
\hline
\rowcolor{lightgray!10} Strategy  & $(p,N)$ & Newton & Avg.~FGMRES & Avg.~GMRES & Avg.~$A^{-1}_\alpha$\\
\hline 
\cref{fig:solver2} &   $(2,4^3)$&   16 & - & 1.69 &1.81  \\
\cref{fig:solver2} &   $(2,8^3)$&   16 & - & 2.75&2.88 \\
\cref{fig:solver2} &   $(6,4^3)$&   16 & - &13.44 &13.56\\
\cref{fig:solver2} &   $(6,8^3)$&   16 & - & 9.5 &9.63 \\
\cref{fig:solver2} &   $(10,4^3)$&   16 & - & 15.13&15.25  \\
\hline
$P_D$ &   $(2,4^3)$&   16 & 3.44  & -& 3.44 \\
$P_D$ &   $(2,8^3)$&   16 & 8.56  & -& 8.56 \\
$P_D$ &   $(6,4^3)$&   16 & 56.88 & - & 56.88  \\
$P_D$  &   $(6,8^3)$&   16 & 69.13  & -& 69.13 \\
$P_D$  &   $(10,4^3)$&   16 & 77.75  & -& 77.75 \\
\hline
$P_LP_D$ &   $(2,4^3)$&   16 & 2.56 & -&  2.56 \\
$P_LP_D$ &   $(2,8^3)$&   16 &  3.94 & -& 3.94 \\
$P_LP_D$ &   $(6,4^3)$&   16 &  30.25 & -&  30.25\\
$P_LP_D$  &   $(6,8^3)$&   16 &  38.88 & -&  38.88 \\
$P_LP_D$  &   $(10,4^3)$&   16 & 43.25  & -&  43.25\\
\hline
$P_F$-(10) &   $(2,4^3)$&   16 & 1.00& 2.36  & 4.36 \\
$P_F$-(10) &   $(2,8^3)$&   16 &  1.00& 5.50& 7.50 \\
$P_F$-(10) &   $(6,4^3)$&   16 &  2.44& 23.63  & 28.50 \\
$P_F$-(10)  &   $(6,8^3)$&   16 &  1.81&  17.31& 20.94 \\
$P_F$-(10)  &   $(10,4^3)$&   16 & 2.81 & 27.31 & 32.94 \\
\hline
\end{tabular}}
\caption{Six iterative solver strategies for the three-dimensional obstacle problem of \Cref{sec:examples:3d}. To invert the Newton linear systems we consider a Schur decomposition as in \cref{fig:solver2}, as well as an outer FGMRES solver block preconditioned as in \cref{fig:solver} with $P=P_D, P_LP_D$ or $P_F$ with $S^{-1}$ either approximated by $\hat{S}^{-1}$ or by inner $\hat{S}$-preconditioned GMRES iterations which are terminated at an absolute or relative tolerance of $10^{-5}$ or after $n$ iterations (denoted by -($n$)). We report the number of Newton iterations as well as the average FGMRES, inner GMRES and $A_\alpha^{-1}$-applications per Newton step.} 
\label{tab:3d}
\end{table}

\section{Conclusions}
\label{sec:conclusions}

In this paper we discretized the latent variable proximal point method \cite{keith2023,dokken2024} with the hierarchical $p$-FEM basis \cite{Babuska1981a,Schwab1998,knook2024quasi} to construct the hierarchical proximal Galerkin (hpG) algorithm. The hpG algorithm is a high-order solver for variational problems with pointwise obstacle- and gradient-type inequality constraints on tensor-product domains. As $p \to \infty$, the choice of FEM basis retains sparsity in the discretized Newton systems and also admits a block preconditioner that, when coupled with an outer FGMRES or inner GMRES iterative solver, bounds the inner Krylov iterations with polylogarithmic growth in the worst case. Moreover, we observed that the outer Newton iterations are experimentally \emph{mesh} and \emph{degree independent}. The method is amenable to $hp$-adaptivity techniques which were explored in a one-dimensional obstacle problem example. We also successively apply the solver to a two- and three-dimensional obstacle problem, a two-dimensional gradient-constrained problem as well as a thermoforming problem, an example of an obstacle-type quasi-variational inequality. We consider discretizations with up to partial degree $p=82$ on each element.

We also compare wall-clock timings of our solver with low-degree discretizations and note that, for the regimes of $h$ and $p$ considered, we found high-order discretizations deliver up to 100 times faster solves to achieve the same errors. This discovery is, perhaps, in contrast to what the community has come to expect.

We now outline some extensions:
\begin{itemize}
\itemsep=-2pt
\item \textbf{Non-Cartesian cells.} The discretization relies on a tensor-product structure of the mesh whenever $d \geq 2$. Non-Cartesian cells are an open problem. For domains that can be built from convex domain building blocks, then equivalent operator preconditioning may offer a viable option \cite{axelsson2009}, \cite[Sec.~2.7]{Brubeck2022}. For domains without structure, further developments are required to construct competitive high-order discretizations.
\item \textbf{$p$-multigrid.} An investigation into $p$-multigrid is worthwhile and may allow one to consider significantly finer meshes \cite{Brubeck2022}.%
\item \textbf{Alternative problems.} Minimizing the Dirichlet energy constitutes the simplest PDE structure one might consider. Proximal Galerkin has been shown to be a flexible framework, as evidenced in \cite{dokken2024}, and can be extended to tackle problems in contact mechanics, fracture,  multiphase species, quasi-variational inequalities, and the Monge--Amp\`ere equation among others. %
\end{itemize}

\section*{Acknowledgements}
I would like to thank Michael Hinterm\"uller, Brendan Keith, Patrick Farrell, Thomas Surowiec, J\o rgen Dokken, Sheehan Olver, Richard M.~Slevinsky, Jishnu Bhattacharya, Lothar Banz, Andreas Schr\"oder, and Pablo Brubeck for discussions and comments that have significantly aided with the writing of this manuscript.

I also acknowledge the funding support by the Deutsche Forschungsgemeinschaft (DFG, German Research Foundation) under Germany's Excellence Strategy -- The Berlin Mathematics Research Center MATH+ (EXC-2046/1, project ID: 390685689).

\printbibliography
\end{document}